\documentclass[10pt]{article}

\usepackage{latexsym}
\usepackage{amsfonts}
\usepackage{amsmath}
\usepackage{mathrsfs}  
\usepackage{dsfont}    
\usepackage{bbold}     
\usepackage[english]{babel}
\usepackage{caption}
\usepackage{epsfig}
\usepackage{float}
\usepackage{subfigure}

\newtheorem{theorem}{Theorem}

\newtheorem{lemma}{Lemma}

\newtheorem{hypothesis}{Hypothesis}

\newcommand{\zerarcounters}{\setcounter{equation}{0}}

\newcommand{\ZZZ}{\mathds{Z}}

\newcommand{\NNN}{\mathds{N}}
\newcommand{\QQQ}{\mathds{Q}}
\newcommand{\RRR}{\mathds{R}}
\newcommand{\TTT}{\mathds{T}}

\newcommand{\FF}{{\mathcal F}}
\newcommand{\GG}{{\mathcal G}}

\newcommand{\MM}{{\mathcal M}}

\newcommand{\gotr}{{\mathfrak r}}

\newcommand{\gotv}{{\mathfrak v}}

\newcommand{\gotw}{{\mathfrak w}}

\newcommand{\gotC}{{\mathfrak C}}

\newcommand{\prova}{\noindent{\it Proof. }}

\newcommand{\eps}{\varepsilon}
\newcommand{\al}{\alpha}
\newcommand{\de}{\delta}

\newcommand{\ka}{\kappa}
\newcommand{\n}{\nu}
\newcommand{\m}{\mu}

\newcommand{\g}{\gamma}
\newcommand{\om}{\omega}

\newcommand{\s}{\sigma}

\newcommand{\Val}{{\rm Val}}

\def\qed{\hfill\raise1pt\hbox{\vrule height5pt width5pt depth0pt}}

\def\ins#1#2#3{\vbox to0pt{\kern-#2 \hbox{\kern#1 #3}\vss}\nointerlineskip}

\begin{document}

\title{\bf Bifurcation curves of subharmonic solutions}
\author{
\bf Guido Gentile$^\dagger$, 
Michele V. Bartuccelli$^\ast$, and Jonathan H.B. Deane$^\ast$
\vspace{2mm}
\\ \small 
$^\dagger$Dipartimento di Matematica, Universit\`a di Roma Tre, Roma,
I-00146, Italy.
\\ \small 
E-mail: gentile@mat.unirom3.it
\\ \small 
$^\ast$Department of Mathematics and Statistics,
University of Surrey, Guildford, GU2 7XH, UK.
\\ \small 
E-mails:
m.bartuccelli@surrey.ac.uk, j.deane@surrey.ac.uk
}

\date{}

\maketitle

\begin{abstract}
We revisit a problem considered by Chow and Hale on the existence of
subharmonic solutions for perturbed systems. In the analytic setting,
under more general (weaker) conditions, we prove their results on the
existence of bifurcation curves from the nonexistence to the existence
of subharmonic solutions. In particular our results apply also when one
has degeneracy to first order --- i.e. when the subharmonic Melnikov
function vanishes identically. Moreover we can deal as well with the
case in which degeneracy persists to arbitrarily high orders, in the
sense that suitable generalisations to higher orders of the
subharmonic Melnikov function are also identically zero. In general the
bifurcation curves are not analytic, and even when they are smooth
they can form cusps at the origin: we say in this case that the curves
are degenerate as the corresponding tangent lines coincide.
The technique we use is completely different from that of Chow and Hale,
and it is essentially based on rigorous perturbation theory.
\end{abstract}

%\newpage

%\tableofcontents

%\newpage

%%%%%%%%%%%%%%%%%%%%%%%%%%%%%%%%%%%%%%%%%%%%%%%%%%%%%%%%%%%%%%%%%%%%%%%%%
%%%%%%%%%%%%%%%%%%%%%%%%%%%%%%%%%%%%%%%%%%%%%%%%%%%%%%%%%%%%%%%%%%%%%%%%%
\zerarcounters
\section{Introduction}\label{sec:1}
%%%%%%%%%%%%%%%%%%%%%%%%%%%%%%%%%%%%%%%%%%%%%%%%%%%%%%%%%%%%%%%%%%%%%%%%%
%%%%%%%%%%%%%%%%%%%%%%%%%%%%%%%%%%%%%%%%%%%%%%%%%%%%%%%%%%%%%%%%%%%%%%%%%

Subharmonic bifurcations have been extensively studied in the
literature \cite{CH,GH}. The problem can be formulated as follows.
Consider a two-dimensional autonomous system, and suppose that
it has a periodic orbit of period $T=2\pi p/q$, where $p,q$
are relatively prime integers. Then one can be interested in studying
whether, under the action of a periodic perturbation of period $2\pi$,
some periodic solutions exist. Solutions with this property
are called subharmonic solutions of order $q/p$.

Assume also that the perturbation depends on two parameters.
A typical situation is when dissipation is present in the
system \cite{HT1,L}; in this case two parameters naturally arise:
the magnitude of the perturbation and the damping coefficient.
An interesting problem is then to study the region in the
space of parameters where subharmonic solutions can occur
and to determine the bifurcation curves, which divide
the regions of existence and non-existence of these solutions.
Such a problem has been considered for instance by Chow and
Hale \cite{CH}. They found that, under suitable assumptions on the
unperturbed system (essentially a local anisochronicity condition)
and on the perturbation, the bifurcation curves exist,
are smooth and intersect with distinct tangent lines at the origin.
The condition on the perturbation, if one takes the magnitude
of the perturbation as one of the parameters, can be formulated in
terms of the so-called subharmonic Melnikov function \cite{M,GH}.
It requires in particular that this function depends explicitly
on the initial phase $t_{0}$ of the unperturbed
periodic solutions which persist under perturbation.

In this paper we recover the same result by Chow and Hale,
in the analytic setting, and we show that the condition
on the perturbation can be weakened. In particular
the subharmonic Melnikov function can be independent of $t_{0}$.
As a consequence the bifurcation curves can be degenerate,
in the sense that they can have the same tangent at the origin,
so that they form a cusp at the origin. Moreover,
in general, they are not smooth. However if some
further assumption is made they turn out to be analytic.

In the case of dissipative systems in the presence of forcing,
such as those studied by Hale and T\'aboas \cite{HT1,CH},
our result is significantly stronger as it requires no
assumption at all on the periodic perturbation.
In particular we find the following result.
Given any one-dimensional anisochronous mechanical system perturbed
by a periodic forcing of magnitude $\eps$ and in the presence of
dissipation, there can be analytic subharmonic solutions of order $q/p$
only if the dissipation coefficient $\g$ is below a threshold
value $\g_{0}(q/p,\eps)$. Here we show that for any rational value
$p/q$ there is an integer exponent $m=m(q/p) \in \RRR^{*}$ such
that $\g_{0}(q/p,\eps)=O(\eps^{m})$. This can be related,
in a more general context, to a conjecture proposed in \cite{BBDGG}.
Moreover the case $m(p/q)=\infty$ corresponds to infinitely many
cancellations, one at each perturbation order, which makes such a case
very unlikely. Therefore, up to these exceptional cases,
we can say that any resonant torus with frequency commensurate
with the frequency of the forcing term admits subharmonic solutions
of the corresponding order. In other words, existence of any
subharmonic solutions holds without making any assumption on the
periodic perturbation, other than smoothness.

Our method is completely different from that of Chow and Hale.
It is based on perturbation theory. More precisely we study
the perturbation series of the subharmonic solutions:
first we find conditions sufficient for these series to
be well-defined to all orders, then we prove that
if the perturbation is small enough convergence of the
series can be proved. Technically, this is achieved by
using the tree formalism, which has been originally introduced
by Gallavotti \cite{G}, inspired by a pioneering paper by
Eliasson \cite{E}, and thereafter has been applied in a long series
of papers on KAM theory \cite{BG0,GG1,G1,G2,GG2,GCB,GBD1,GBD2};
see also \cite{GBG} for a review. We note that with respect
to these papers in our case the analysis is much easier
as we deal with periodic solutions instead of quasi-periodic
solutions. In this respect our analysis could
be considered as a propaedeutic introduction to the
tree formalism, in a case in which there is no small divisors
problem, so that no multiscale analysis has to be introduced;
see also \cite{BG1,BG2} for a similar situation.
In particular Chow and Hale's assumptions on the perturbation
reflect a case in which a first order analysis is enough
to deduce existence of subharmonic solutions.
By contrast our results allow the analysis of cases in
which it can be necessary to go beyond the first order,
in principle to arbitrarily high orders.

We also argue that in physical applications it can be essential
to have a stronger result. Indeed, in a concrete example in which,
for instance, the perturbation is a trigonometric polynomial,
Chow and Hale's assumptions on the perturbation, even if they
are generic, fail to be satisfied for most values of the periods $T$.
For those values a first order condition is not sufficient
to detect the existence of the subharmonic solution,
and one must go to higher orders. The numerical simulations
performed in \cite{BBDGG} for a driven quartic oscillator
in the presence of dissipation  show that this is necessary
if one wants to explain the numerical findings for some values
of the parameters.

Our results should be compared also with \cite{C,CL}, where a different
scenario, such as the persistence of the whole invariant manifold
corresponding to the resonant torus, arises in a case
in which the subharmonic Melnikov function vanishes identically.
Our analysis shows that a situation of this kind is highly non-generic.

%%%%%%%%%%%%%%%%%%%%%%%%%%%%%%%%%%%%%%%%%%%%%%%%%%%%%%%%%%%%%%%%%%%%%%%%%
%%%%%%%%%%%%%%%%%%%%%%%%%%%%%%%%%%%%%%%%%%%%%%%%%%%%%%%%%%%%%%%%%%%%%%%%%
\zerarcounters
\section{Main results}\label{sec:2}
%%%%%%%%%%%%%%%%%%%%%%%%%%%%%%%%%%%%%%%%%%%%%%%%%%%%%%%%%%%%%%%%%%%%%%%%%
%%%%%%%%%%%%%%%%%%%%%%%%%%%%%%%%%%%%%%%%%%%%%%%%%%%%%%%%%%%%%%%%%%%%%%%%%

Consider the ordinary differential equation
\begin{equation}
\begin{cases}
\dot \al = \om(A) + \eps F(\al,A,C,t) , \\
\dot A = \eps G(\al,A,C,t) , \end{cases}
\label{eq:2.1} \end{equation}
where $(\al,A) \in \MM := \TTT \times W$, with $W \subset\RRR$ an
open set, the map $A\to \om(A)$ is analytic in $A$,
and the functions $F$ and $G$ depend analytically
on their arguments and are $2\pi$-periodic in $\al$ and $t$.
Finally, $\eps$, $C$ are two real parameters.

One could also introduce a further (analytic) dependence on $\eps$
in the functions $F$ and $G$, and the forthcoming analysis
could be easily performed with some trivial adaptations.
Therefore all the results and theorems stated below and in the
next sections hold unchanged in that case too. Then, the
formulation given in \cite{CH} is recovered, as a particular case,
by introducing the parameter $\g=\eps C$, and setting
$\m=(\m_{1},\m_{2})$, with $\m_{1}=\eps$ and $\m_{2}=\g$.

For $\eps=0$ the variable $A$ is kept fixed at some value $A_{0}$,
while $\al$ rotates with constant angular velocity $\om(A_{0})$.
Hence the motion of the variables $(\al,A,t)$ is quasi-periodic,
and reduces to a periodic motion whenever $\om(A_{0})$ becomes
commensurate with $1$. Define $\al_{0}(t)=\om(A_{0})t$ and $A_{0}(t)=
A_{0}$: in the \textit{extended phase-space} $\MM \times \RRR$ the
solution $(\al_{0}(t),A_{0}(t),t+t_{0})$ describes an invariant torus,
which is uniquely determined by the ``energy'' $A_{0}$. If $\om(A_{0})$
is rational we say that the torus is \textit{resonant}.
The parameter $t_{0}$ will be called the \textit{initial phase}:
it fixes the initial datum on the torus.

As a particular case we can consider that $(A,\al)$ are canonical
coordinates (action-angle coordinates), but the formulation we are
giving here is more general. In general all non-resonant tori
are completely destroyed under perturbation, if no further hypotheses
are made on the perturbations $F,G$ (such as that the
full system is Hamiltonian). Also the resonant tori disappear,
but some remnants are left: indeed usually a finite number
of periodic orbits, called \textit{subharmonic solutions},
lying on the unperturbed torus, can survive under perturbation.

Denote by $T_{0}(A)=2\pi/\om(A)$ the period
of the trajectories on an unperturbed torus, and define
$\om'(A):={\rm d}\om(A)/{\rm d}A$. If $\om(A_{0})=p/q \in
\QQQ$, call $T=T(A_{0})=2\pi q$ the period of the trajectories
in the extended phase space. We shall call $q/p$ the \textit{order}
of the corresponding subharmonic solutions. Define
\begin{equation}
M(t_{0},C) = \frac{1}{T} \int_{0}^{T} {\rm d} t \,
G(\al_{0}(t),A_{0},C,t+t_{0}) ,
\label{eq:2.2} \end{equation}
which is called the \textit{subharmonic Melnikov function}.
Here and in the following we do not write explicitly the dependence
of the subharmonic Melnikov function on $A_{0}$,
which is fixed once and for all. Note that $M(t_{0},C)$
is $2\pi$-periodic in $t_{0}$.

We make the following assumptions on the resonant torus
with energy $A_{0}$.

%%%%%%%%%%%%%%%%%%%%%%%%%%%%%%%%%%%%%%%%%%%%%%%%%%%%%%%%%%%%%%%%%%%%%%%%%
\begin{hypothesis} \label{hyp:1}
One has $\om'(A_{0}) \neq 0$.
\end{hypothesis}
%%%%%%%%%%%%%%%%%%%%%%%%%%%%%%%%%%%%%%%%%%%%%%%%%%%%%%%%%%%%%%%%%%%%%%%%%

%%%%%%%%%%%%%%%%%%%%%%%%%%%%%%%%%%%%%%%%%%%%%%%%%%%%%%%%%%%%%%%%%%%%%%%%%
\begin{hypothesis} \label{hyp:2}
There exists an analytic curve $t \to C_{0}(t)$ from $[0,2\pi)$
to $\RRR$ such that $M(t_{0},C_{0}(t_{0}))=0$ and
$\partial M(t_{0},C_{0}(t_{0}))/\partial C \neq 0$
for all $t_{0}\in[0,2\pi)$.
\end{hypothesis}
%%%%%%%%%%%%%%%%%%%%%%%%%%%%%%%%%%%%%%%%%%%%%%%%%%%%%%%%%%%%%%%%%%%%%%%%%

The function $C_{0}(t_{0})$ is also $2\pi$-periodic in $t_{0}$.
We prove the following result. We prefer to state the result in
terms of the parameter $\g=\eps C$ --- instead of $C$ --- to make more
transparent the relation with \cite{CH}.

%%%%%%%%%%%%%%%%%%%%%%%%%%%%%%%%%%%%%%%%%%%%%%%%%%%%%%%%%%%%%%%%%%%%%%%%%
\begin{theorem} \label{thm:1}
Consider the system (\ref{eq:2.1}) and assume that Hypotheses
\ref{hyp:1} and \ref{hyp:2} hold for the resonant torus with energy
$A_{0}$ such that $\om(A_{0})=p/q$. There exist $\eps_{0} > 0$ and two
continuous functions $\g_{1}(\eps)$ and $\g_{2}(\eps)$, with
$\g_{1}(0)=\g_{2}(0)$, $\g_{1}(\eps)\ge \g_{2}(\eps)$ for $\eps\ge 0$
and $\g_{1}(\eps)\le \g_{2}(\eps)$ for $\eps\le 0$,
such that (\ref{eq:2.1}) has at least one
subharmonic solution of order $q/p$ for $\g_{2}(\eps) \le \eps C \le
\g_{1}(\eps)$ when $\eps\in(0,\eps_{0})$ and for $\g_{1}(\eps) \le
\eps C \le \g_{2}(\eps)$ when $\eps\in(-\eps_{0},0)$.
\end{theorem}
%%%%%%%%%%%%%%%%%%%%%%%%%%%%%%%%%%%%%%%%%%%%%%%%%%%%%%%%%%%%%%%%%%%%%%%%%

The situation is depicted in Figure 1, in a case in which the two
functions $\g_{1}$ and $\g_{2}$ are smooth. The graphs described
by the two functions are called the \textit{bifurcation curves} of
the subharmonic solutions: they divide the plane into two disjoint sets
such that only in one of them there are analytic subharmonic solutions.

%%%%%%%%%%%%%%%%%%%%%%%%%%%%%%%%%%%%%%%%%%%%%%%%%%%%%%%%%%%%%%%%%%%%%%%%%
%FIGURE 1
%%%%%%%%%%%%%%%%%%%%%%%%%%%%%%%%%%%%%%%%%%%%%%%%%%%%%%%%%%%%%%%%%%%%%%%%%
\begin{figure}[htbp]
\begin{centering}
\ins{358pt}{-094pt}{$\eps$}
\ins{276pt}{-024pt}{$\g_{1}(\eps)$}
\ins{270pt}{-128pt}{$\g_{2}(\eps)$}
\ins{170pt}{-040pt}{$\g_{2}(\eps)$}
\ins{170pt}{-124pt}{$\g_{1}(\eps)$}
\ins{214pt}{-004pt}{$\g$}
\includegraphics*[angle=0,width=4in]{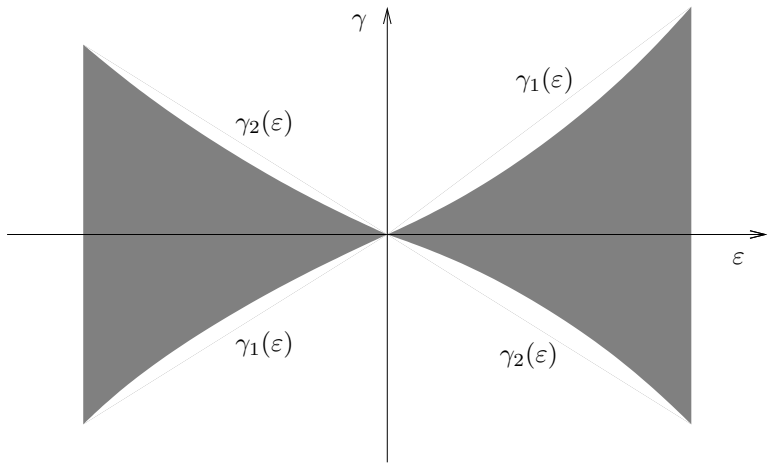}
\caption{\label{fig2-1}
Set of existence (grey region) of subharmonic solutions
in the plane $(\eps,\g)$, in a case in which
the two bifurcation curves $\eps \to \g_{1}(\eps)$ and
$\eps \to \g_{2}(\eps)$ are smooth and
have different tangent lines at the origin.}
\end{centering}
\end{figure}
%%%%%%%%%%%%%%%%%%%%%%%%%%%%%%%%%%%%%%%%%%%%%%%%%%%%%%%%%%%%%%%%%%%%%%%%%

In general the functions $\g_{1}$ and $\g_{2}$ are not smooth.
However, if some further assumptions are made on the subharmonic
Melnikov function, smoothness (in fact analyticity) in $\eps$ can be
obtained. Denote by $C_{0}'(t_{0})$ and $C_{0}''(t_{0})$ the first and
second derivatives of the function $C_{0}(t_{0})$ with respect to $t_{0}$.

%%%%%%%%%%%%%%%%%%%%%%%%%%%%%%%%%%%%%%%%%%%%%%%%%%%%%%%%%%%%%%%%%%%%%%%%%
\begin{hypothesis} \label{hyp:3}
If $t_{m}$ and $t_{M}$ are the values
in $[0,2\pi)$ for which the function $C_{0}(t_{0})$ attains
its minimum and its maximum, respectively,
then $C_{0}''(t_{m}) C_{0}''(t_{M}) \neq 0$.
\end{hypothesis}
%%%%%%%%%%%%%%%%%%%%%%%%%%%%%%%%%%%%%%%%%%%%%%%%%%%%%%%%%%%%%%%%%%%%%%%%%

The following result holds.

%%%%%%%%%%%%%%%%%%%%%%%%%%%%%%%%%%%%%%%%%%%%%%%%%%%%%%%%%%%%%%%%%%%%%%%%%
\begin{theorem} \label{thm:2}
Consider the system (\ref{eq:2.1}) and assume that Hypotheses
\ref{hyp:1}, \ref{hyp:2} and \ref{hyp:3} hold for the resonant torus
with energy $A_{0}$ such that $\om(A_{0})=p/q$. There exist
$\eps_{0} > 0$ and two functions $\g_{1}(\eps)$ and
$\g_{2}(\eps)$, analytic for $|\eps|<\eps_{0}$, with
$\g_{1}(0)=\g_{2}(0)$, $\g_{1}(\eps) > \g_{2}(\eps)$ for $\eps > 0$
and $\g_{1}(\eps) < \g_{2}(\eps)$ for $\eps < 0$,
and with different tangent lines at the origin, such that (\ref{eq:2.1})
has at least one subharmonic solution of order $q/p$
for $\g_{2}(\eps) \le \eps C \le \g_{1}(\eps)$ when $\eps\in(0,\eps_{0})$
and for $\g_{1}(\eps) \le \eps C \le \g_{2}(\eps)$
when $\eps\in(-\eps_{0},0)$.
\end{theorem}
%%%%%%%%%%%%%%%%%%%%%%%%%%%%%%%%%%%%%%%%%%%%%%%%%%%%%%%%%%%%%%%%%%%%%%%%%

Theorem \ref{thm:2} is analogous to Theorem 2.1 of \cite{CH},
Section 11 --- in the analytic setting instead of the differentiable
one --- while Theorem \ref{thm:1} requires fewer hypotheses.
In particular it applies when Chow and Hale's $h_{k}(\al)$
function vanishes identically. In that case the graphs of the two
functions $\g_{1}$ and $\g_{2}$ form a cusp at the origin: we refer to this
situation as a case of \textit{degenerate bifurcation curves}, see
Figure 2. We shall also see in Section \ref{sec:4} that in fact,
under weaker assumptions than those made in Hypothesis \ref{hyp:3},
we can find smoothness of the bifurcation curves, in the following
sense: under suitable assumptions there exist two analytic
functions $\widetilde \g_{1}(\eps)$ and $\widetilde \g_{2}(\eps)$
such that $\g_{1}(\eps)= \max\{\widetilde\g_{1}(\eps),
\widetilde \g_{2}(\eps)\}$ and $\g_{2}(\eps)=\min\{
\widetilde\g_{1}(\eps),\widetilde \g_{2}(\eps)\}$
for $\eps>0$, and $\g_{1}(\eps)= \min\{\widetilde\g_{1}(\eps),
\widetilde \g_{2}(\eps)\}$ and $\g_{2}(\eps)= \max\{
\widetilde\g_{1}(\eps),\widetilde \g_{2}(\eps)\}$ for $\eps<0$.
We refer to Hypothesis \ref{hyp:4} and Theorem \ref{thm:3}
in Section \ref{sec:4} for a precise formulation of the results.

%%%%%%%%%%%%%%%%%%%%%%%%%%%%%%%%%%%%%%%%%%%%%%%%%%%%%%%%%%%%%%%%%%%%%%%%%
%FIGURE 2
%%%%%%%%%%%%%%%%%%%%%%%%%%%%%%%%%%%%%%%%%%%%%%%%%%%%%%%%%%%%%%%%%%%%%%%%%
\begin{figure}[htbp]
\begin{centering}
\ins{358pt}{-094pt}{$\eps$}
\ins{310pt}{-000pt}{$\g_{1}(\eps)$}
\ins{300pt}{-144pt}{$\g_{2}(\eps)$}
\ins{120pt}{-004pt}{$\g_{2}(\eps)$}
\ins{136pt}{-144pt}{$\g_{1}(\eps)$}
\ins{214pt}{-004pt}{$\g$}
\includegraphics*[angle=0,width=4in]{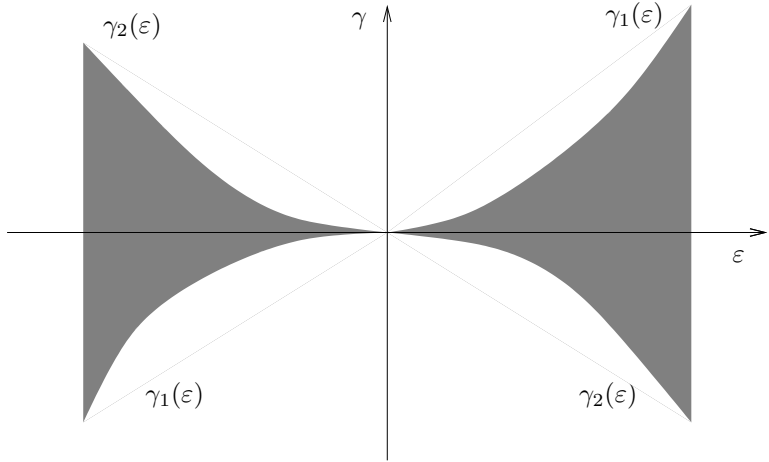}
\caption{\label{fig2-2}
In general the bifurcation curves are not smooth. Even when they
are smooth they can have the same tangent at the origin:
in this case we say that the bifurcation curves are degenerate.
The grey region in the figure represents a case in which the
bifurcation curves are smooth and both of them
have tangent lines parallel to the $\eps$-axis.}
\end{centering}
\end{figure}
%%%%%%%%%%%%%%%%%%%%%%%%%%%%%%%%%%%%%%%%%%%%%%%%%%%%%%%%%%%%%%%%%%%%%%%%%

We shall see in Section \ref{sec:4} --- cf. Theorem \ref{thm:4} --- that
for $p=1$ one has at least $2q$ subharmonic solutions of order $q$
as far as $\min\{\g_{1}(\eps),\g_{2}(\eps)\} < \g < \max\{\g_{1}(\eps),
\g_{2}(\eps)\}$ and at least $q$ subharmonic solutions of order $q$
when $(\eps,\g)$ belongs to one of the bifurcation curves,
that is when either $\g=\g_{1}(\eps)$ or $\g=\g_{2}(\eps)$. This
agrees with Chow and Hale's Theorem 2.1 in \cite{CH} in the cases
in which it applies.

Possible extensions of Chow and Hale's results could be looked
for in another direction, such as that of relaxing the hypothesis
on the unperturbed system. This problem has been studied,
for instance, in \cite{HT2,HT3,RC}.

The rest of the paper is organised as follows.
Sections~\ref{sec:3} and \ref{sec:4} are devoted to the
proof of Theorems \ref{thm:1} and \ref{thm:2}.
More precisely, in Section \ref{sec:2} we show the existence
of a subharmonic solution in the form of a formal power series,
while in Section \ref{sec:4} we prove the convergence of the series,
and we also state Theorem \ref{thm:3}, which provides the
aforementioned extensions of Theorem \ref{thm:2}, and Theorem
\ref{thm:4} on the minimal number of subharmonic solutions of order $q$.
In Section \ref{sec:5} we discuss, as an application of our results,
the case of a forced one-dimensional system in the presence of
dissipation: this will lead to Theorems \ref{thm:5} and \ref{thm:6}
which extend the results of Hale and T\'aboas \cite{HT1}.

Some final comments, and a comparison with the standard Melnikov
theory \cite{M,GH}, are provided in Section \ref{sec:6}. In particular
we shall formulate two theorems: Theorem \ref{thm:7} corresponds
to the result usually discussed in the literature \cite{GH}, while
Theorem \ref{thm:8} provides an extension to degenerate situations.

%%%%%%%%%%%%%%%%%%%%%%%%%%%%%%%%%%%%%%%%%%%%%%%%%%%%%%%%%%%%%%%%%%%%%%%%%
%%%%%%%%%%%%%%%%%%%%%%%%%%%%%%%%%%%%%%%%%%%%%%%%%%%%%%%%%%%%%%%%%%%%%%%%%
\zerarcounters
\section{Existence of formal power series for the subharmonic solutions}
\label{sec:3}
%%%%%%%%%%%%%%%%%%%%%%%%%%%%%%%%%%%%%%%%%%%%%%%%%%%%%%%%%%%%%%%%%%%%%%%%%
%%%%%%%%%%%%%%%%%%%%%%%%%%%%%%%%%%%%%%%%%%%%%%%%%%%%%%%%%%%%%%%%%%%%%%%%%

We look for subharmonic solutions of (\ref{eq:2.1})
which are analytic in $\eps$. First, we shall try
to find solutions in the form of formal power series in $\eps$
\begin{equation}
\al(t) = \al(t,\eps) = \sum_{k=0}^{\infty}
\eps^{k} \al^{(k)}(t) , \qquad
A(t) = A(t;\eps) = \sum_{k=0}^{\infty}
\eps^{k} A^{(k)}(t) ,
\label{eq:3.1} \end{equation}
where $\al^{(0)}(t)=\om(A_{0}) \, t$ and $A^{(0)}(t)=A_{0}$,
with $\om(A_{0})=p/q$, and the functions $\al^{(k)}(t)$
and $A^{(k)}(t)$, periodic with period $T=2\pi p$
for all $k\in\NNN$, are to be determined. We shall see that this
will be possible provided the parameter $C$ is chosen
as a function of $\eps$, again in the form
of a formal power series in $\eps$
\begin{equation}
C = C(\eps) = \sum_{k=0}^{\infty} \eps^{k} C^{(k)} .
\label{eq:3.2} \end{equation}
Moreover both the solution $(\al(t),A(t))$ and the constant $C$
will be found to depend on the initial phase $t_{0}$: in particular
one has $C(\eps)=C(\eps,t_{0})$ such that $C(\eps,t_{0}+2\pi)=
C(\eps,t_{0})$ and $C^{(0)}=C_{0}(t_{0})$, and, as we shall see,
a sufficient condition for formal solvability to hold is
that Hypotheses \ref{hyp:1} and \ref{hyp:2} are satisfied.

If we introduce the decompositions (\ref{eq:3.1}) and (\ref{eq:3.2})
into (\ref{eq:2.1}) and we denote with $W(t)$ the Wronskian
matrix for the unperturbed linearised system,
we obtain (cf. \cite{BBDGG} for similar computations)
\begin{equation}
\left( \begin{matrix} \al^{(k)}(t) \\
A^{(k)}(t) \end{matrix} \right)
= W(t) \left( \begin{matrix} \Bar \al^{(k)} \\
\Bar A^{(k)} \end{matrix} \right) +
W(t) \int_{0}^{t} {\rm d} \tau \, W^{-1}(\tau)
\left( \begin{matrix} F^{(k-1)}(\tau) \\
G^{(k-1)}(\tau) \end{matrix} \right) ,
\label{eq:3.3} \end{equation}
where $(\Bar \al^{(k)},\Bar A^{(k)})$ are corrections
to the initial conditions, and
\begin{eqnarray}
& & F^{(k)}(t) = \left[ F(\al,A,C,t+t_{0}) \right]^{(k)} :=
\sum_{m=0}^{\infty} \sum_{\substack{r_{1},r_{2},r_{3} \in \ZZZ_{+}
\\ r_{1}+r_{2}+r_{3}=m}}
\frac{\partial_{1}^{r_{1}} \partial_{2}^{r_{2}}
\partial_{3}^{r_{3}}}{r_{1}!r_{2}!r_{3}!}
F(\al_{0}(t),A_{0},C_{0},t+t_{0}) \nonumber \\
& & \sum_{k_{1}+\ldots+k_{m}=k}
\al^{(k_{1})}(t)\ldots \al^{(k_{r_{1}})}(t) \,
A^{(k_{r_{1}+1})}(t)\ldots A^{(k_{r_{1}+r_{2}})}(t) \,
C^{(k_{r_{1}+r_{2}+1})}\ldots C^{(k_{m})} ,
\label{eq:3.4} \end{eqnarray}
with an analogous definition holding for $G^{(k)}(t)$. Here
and henceforth, given a function of several arguments
we are denoting by $\partial_{k}$ the derivative with respect
to the $k$-th argument; hence, given the function
$F(\al,A,C,t+t_{0})$ we have $\partial_{1}F=\partial F/\partial\al$,
$\partial_{2}F=\partial F/\partial A$,
and $\partial_{3}F=\partial F/\partial C$.
Note that by construction both $F^{(k)}(t)$ and $G^{(k)}(t)$
depend only on the coefficients $\al^{(k')}(t)$,
$A^{(k')}(t)$ and $C^{(k')}$ with $k' \le k$.

The Wronskian matrix appearing in (\ref{eq:3.3}) can be written as
\begin{equation}
W(t) = \left( \begin{matrix}
1 & \om'(A_{0}) t \\ 0 & 1 \end{matrix} \right) .
\label{eq:3.5} \end{equation}
By using (\ref{eq:3.5}) in (\ref{eq:3.3}) we have
\begin{equation}
\begin{cases}
{\displaystyle
\al^{(k)}(t) = \Bar \al^{(k)} + t \, \om'(A_{0})\,\Bar A^{(k)} +
\int_{0}^{t} {\rm d}\tau \, F^{(k-1)}(\tau) + \om'(A_{0})
\int^{t}_{0} {\rm d} \tau \int_{0}^{\tau} {\rm d}\tau'
G^{(k-1)}(\tau') , } \\
{\displaystyle A^{(k)}(t) = \Bar A^{(k)} +
\int_{0}^{t} {\rm d}\tau \, G^{(k-1)}(\tau) . }
\end{cases}
\label{eq:3.6} \end{equation}
We obtain a periodic solution of period $T$ if,
to any order $k\in\NNN$, one has
\begin{equation}
\langle G^{(k-1)} \rangle : = \frac{1}{T}
\int_{0}^{T} {\rm d}\tau \, G^{(k-1)}(\tau) = 0
\label{eq:3.7} \end{equation}
and
\begin{equation}
\om'(A_{0})\,\Bar A^{(k)} + \langle F^{(k-1)} \rangle +
\om'(A_{0}) \langle \GG^{(k-1)} \rangle = 0 , \qquad
\GG^{(k-1)}(t) = \int_{0}^{t} {\rm d} \tau \, G^{(k-1)}(\tau) ,
\label{eq:3.8} \end{equation}
where, given any $T$-periodic function $H$
we denote by $\langle H \rangle$ its mean, as done in (\ref{eq:3.7}).

The parameters $\Bar\al^{(k)}$ are left undetermined,
and we can fix them arbitrarily, as we have the initial phase $t_{0}$
which is still a free parameter. For instance we can
set $\Bar\al^{(k)}=0$ for all $k\in\NNN$ or else
we can define $\Bar\al^{(k)}=\al_{k}(t_{0})$ for $k\in\NNN$,
with the constants $\al_{k}(t_{0})$ to be fixed in the
way which turns out to be more convenient for computations:
we shall see in the next section a reasonable choice.

Therefore, if equation (\ref{eq:3.7}) is satisfied, we have
\begin{equation}
\begin{cases}
{\displaystyle \al^{(k)}(t) = \int_{0}^{t} {\rm d}\tau
\left( F^{(k-1)}(\tau) - \langle F^{(k-1)} \rangle \right) +
\om'(A_{0}) \int^{t}_{0} {\rm d} \tau \left( \GG^{(k-1)}(\tau) -
\langle \GG^{(k-1)} \rangle \right) , } \\
\\
{\displaystyle A^{(k)}(t) = \Bar A^{(k)} + \GG^{(k-1)}(t) , }
\end{cases}
\label{eq:3.9} \end{equation}
with
\begin{equation}
\Bar A^{(k)} = - \frac{\langle F^{(k-1)} \rangle}{\om'(A_{0})} -
\langle \GG^{(k-1)} \rangle = 0 ,
\label{eq:3.10} \end{equation}
which is well-defined as $\om'(A_{0}) \neq 0$ by Hypothesis \ref{hyp:1}.

So, in order to prove the formal solvability of (\ref{eq:2.1})
we have to check whether it is possible to fix the parameter $C$,
as a function of $\eps$ and $t_{0}$, in such a way that
(\ref{eq:3.7}) follows for all $k\ge 1$.

For $k=1$ the condition (\ref{eq:3.7}) reads
\begin{equation}
\langle G^{(0)} \rangle = M(t_{0},C) = 0 ,
\label{eq:3.11} \end{equation}
and we can choose $C=C_{0}(t_{0})$ so that this holds:
this is assured by Hypothesis \ref{hyp:2}.

To higher order $k\ge 1$ we can write
\begin{equation}
G^{(k)}(\al(t),A(t),C,t+t_{0}) =
\partial_{3} G(\al_{0}(t),A_{0},C_{0},t+t_{0}) \, C^{(k)}
+ \Gamma^{(k)}(\al(t),A(t),C,t+t_{0}) ,
\label{eq:3.12} \end{equation}
where the function $\Gamma^{(k)}(\al(t),A(t),C,t+t_{0})$
depends on the coefficients $C^{(k')}$ of $C$ with $k'<k$
(and on the functions $\al^{(k')}(t)$ and $A^{(k')}(t)$ with
$k'\le k$, of course). In other words, in (\ref{eq:3.12})
we have extracted explicitly the only term depending on $C^{(k)}$.
Moreover one has
\begin{equation}
\langle \partial_{3} G(\al_{0}(\cdot),A_{0},C_{0},\cdot+t_{0})
\rangle = \frac{1}{T} \int_{0}^{T} {\rm d}t \, \partial_{3}
G(\al_{0}(t),A_{0},C_{0},t+t_{0}) =
\frac{\partial}{\partial C} M(t_{0},C_{0}),
\label{eq:3.13} \end{equation}
and by Hypothesis \ref{hyp:2} one has
$D(t_{0}) := \partial M(t_{0},C_{0}(t_{0})) / \partial C \neq 0$,
so that (\ref{eq:3.7}) is satisfied provided $C^{(k)}$ is chosen as
\begin{equation}
C^{(k)} = - \frac{1}{D(t_{0})} \langle
\Gamma^{(k)}(\al(\cdot),A(\cdot),C,\cdot+t_{0}) \rangle
\equiv C_{k}(t_{0}) .
\label{eq:3.14} \end{equation}
Therefore we conclude that if we set $C_{0}=C_{0}(t_{0})$
and, for all $k\ge 1$, we choose $\Bar\al^{(k)}=\al_{k}(t_{0})$,
$\Bar A^{(k)}$ according to (\ref{eq:3.10}) and
$C^{(k)}=C_{k}(t_{0})$ according to (\ref{eq:3.14}), we obtain
that in the expansions (\ref{eq:3.1}) the coefficients 
$\al^{(k)}(t)$ and $A^{(k)}(t)$ are well-defined periodic
functions of period $T$. Of course this does not settle
the problem of convergence of the series (\ref{eq:3.1})
and (\ref{eq:3.2}). This will be discussed in the next Section.

%%%%%%%%%%%%%%%%%%%%%%%%%%%%%%%%%%%%%%%%%%%%%%%%%%%%%%%%%%%%%%%%%%%%%%%%%
%%%%%%%%%%%%%%%%%%%%%%%%%%%%%%%%%%%%%%%%%%%%%%%%%%%%%%%%%%%%%%%%%%%%%%%%%
\zerarcounters
\section{Convergence of the series for the subharmonic solutions}
\label{sec:4}
%%%%%%%%%%%%%%%%%%%%%%%%%%%%%%%%%%%%%%%%%%%%%%%%%%%%%%%%%%%%%%%%%%%%%%%%%
%%%%%%%%%%%%%%%%%%%%%%%%%%%%%%%%%%%%%%%%%%%%%%%%%%%%%%%%%%%%%%%%%%%%%%%%%

Here we shall prove that the formal power series found
in Section \ref{sec:3} converge for $\eps$ small enough,
say for $|\eps|<\eps_{0}$ for some $\eps_{0}>0$.
Then for fixed $\eps\in(-\eps_{0},\eps_{0})$ we shall find the range
allowed for $C$ by computing the supremum and the infimum, for
$t_{0}\in[0,2\pi)$ of the function $t_{0} \to C(\eps,t_{0})$.
The bifurcation curves will be defined in terms of the function
$C(\eps,t_{0})$ --- cf. (\ref{eq:3.2}) --- as
\begin{equation}
\g_{1}(\eps) = \eps \sup_{t_{0}\in[0,2\pi)} C(\eps,t_{0}) , \qquad
\g_{2}(\eps) = \eps \inf_{t_{0}\in[0,2\pi)} C(\eps,t_{0}) .
\label{eq:4.1} \end{equation}
In general the functions (\ref{eq:4.1}) are not smooth in $\eps$.
We shall return to this at the end of the section.

To prove convergence of the series (\ref{eq:3.1}) and (\ref{eq:3.2})
it is more convenient to work in Fourier space.
First of all let us define $\om=2\pi/T=1/q$
(note that $\om \neq \om(A_{0})$) and expand
\begin{equation}
F(\al,A,C,t+t_{0}) =
\sum_{\n,\s\in\ZZZ} {\rm e}^{i\n\al} {\rm e}^{i\s(t+t_{0})}
F_{\n,\s}(A,C) ,
\label{eq:4.2} \end{equation}
so that we can write
\begin{equation}
\partial_{1}^{r_{1}} \partial_{2}^{r_{2}} \partial_{3}^{r_{3}}
F(\om(A_{0})\,t,A_{0},C_{0}(t_{0}),t+t_{0}) =
\sum_{\n\in\ZZZ} {\rm e}^{i\n\om t}
\sum_{\substack{\n_{0},\s_{0}\in\ZZZ \\ \n_{0}p+\s_{0}q=\n}}
{\rm e}^{i\s t_{0}} \left( i\n_{0} \right)^{r_{1}}
\partial_{2}^{r_{2}} \partial_{3}^{r_{3}}
F_{\n_{0},\s_{0}}(A_{0},C_{0}(t_{0})) ,
\label{eq:4.3} \end{equation}
and an analogous expression can be obtained with the function $G$
replacing $F$. By the analyticity assumption on the functions
$F$ and $G$, we have the bounds
\begin{eqnarray}
& & \left| \frac{\partial_{2}^{r_{2}} \partial_{3}^{r_{3}}}{r_{2}!r_{3}!}
F_{\n_{0},\s_{0}}(A_{0},C_{0}(t_{0})) \right| \le
P Q_{1}^{r_{1}}Q_{2}^{r_{2}} {\rm e}^{-\ka(|\n_{0}|+|\s_{0}|)} ,
\nonumber \\
& & \left| \frac{\partial_{2}^{r_{2}} \partial_{3}^{r_{3}}}{r_{2}!r_{3}!}
G_{\n_{0},\s_{0}}(A_{0},C_{0}(t_{0})) \right| \le
P Q_{1}^{r_{1}}Q_{2}^{r_{2}} {\rm e}^{-\ka(|\n_{0}|+|\s_{0}|)} ,
\label{eq:4.4} \end{eqnarray}
for suitable positive constants $P,Q_{1},Q_{2},\ka$.

Then, let us write in (\ref{eq:3.1})
\begin{equation}
\al^{(k)}(t) = \sum_{\nu\in\ZZZ} {\rm e}^{i\n\om t} \al^{(k)}_{\n} ,
\qquad A^{(k)}(t) = \sum_{\nu\in\ZZZ} {\rm e}^{i\n\om t} A^{(k)}_{\n} ,
\label{eq:4.5} \end{equation}
so that (\ref{eq:3.9}) becomes
\begin{equation}
\al^{(k)}_{\n} = \frac{F^{(k-1)}_{\n}}{i\om\nu} +
\om'(A_{0}) \frac{G^{(k-1)}_{\n}}{(i\om\nu)^{2}} , \qquad
A^{(k)}_{\n} = \frac{G^{(k-1)}_{\n}}{i\om\nu} ,
\label{eq:4.6} \end{equation}
for all $\n\neq0$, whereas for $\n=0$ one has
\begin{equation}
\al^{(k)}_{0} = \al_{k}(t_{0}) - \sum_{\substack{\n\in\ZZZ \\ \n\neq 0}}
\frac{F^{(k-1)}_{\n}}{i\om\nu} - \om'(A_{0})
\sum_{\substack{\n\in\ZZZ \\ \n\neq 0}}
\frac{G^{(k-1)}_{\n}}{(i\om\nu)^{2}} , \qquad
A^{(k)}_{0} = \Bar A^{(k)} - \sum_{\substack{\n\in\ZZZ \\ \n\neq 0}}
\frac{G^{(k-1)}_{\n}}{i\om\nu} = -\frac{F^{(k-1)}_{0}}{\om'(A_{0})} ,
\label{eq:4.7} \end{equation}
with $\al_{k}(t_{0})$ so far arbitrary and $\Bar A^{(k)}$
given by (\ref{eq:3.10}). The Fourier coefficients $F^{(k-1)}_{\n}$
and $G^{(k-1)}_{\n}$ can be read from (\ref{eq:3.4}) and the analogous
expression for $G^{(k)}(t)$. Hence one has
\begin{eqnarray}
& & \null \hskip-1.4truecm
F^{(k)}_{\n} = \left[ F(\al,A,C,t+t_{0}) \right]^{(k)}_{\n} =
\sum_{m=0}^{\infty} \sum_{\substack{r_{1},r_{2},r_{3} \in \ZZZ_{+} \\
r_{1}+r_{2}+r_{3}=m}} \hskip.1truecm \sum_{\substack{
\n_{0},\s_{0},\n_{1},\ldots,\n_{r_{1}+r_{2}} \in\ZZZ \\
\n_{0}p+\s_{0}q+\n_{1}+\ldots+\n_{r_{1}+r_{2}}=\n}}
\frac{(i\n_{0})^{r_{1}}}{r_{1}!} 
{\rm e}^{i\s_{0} t_{0}} \nonumber \\
& & \null \hskip-1.4truecm
\qquad \frac{\partial_{2}^{r_{2}} \partial_{3}^{r_{3}}}{r_{2}!r_{3}!}
 F_{\n_{0}}(A_{0},C_{0}(t_{0})) \sum_{k_{1}+\ldots+k_{m}=k}
\al^{(k_{1})}_{\n_{1}}\ldots \al^{(k_{r_{1}})}_{\n_{r_{1}}} \,
A^{(k_{r_{1}+1})}_{\n_{r_{1}+1}}
\ldots A^{(k_{r_{1}+r_{2}})}_{\n_{r_{1}+r_{2}}} \,
C^{(k_{r_{1}+r_{2}+1})}\ldots C^{(k_{m})} ,
\label{eq:4.8} \end{eqnarray}
and an analogous definition holds for $G^{(k)}_{\n}$.

Furthermore one has
\begin{equation}
C^{(k)} = -\frac{1}{D(t_{0})} \Gamma^{(k)}_{0}
\label{eq:4.9} \end{equation}
where
\begin{eqnarray}
& & \null \hskip-1.4truecm
\Gamma^{(k)}_{0} = \left[ \Gamma(\al,A,C,t+t_{0}) \right]^{(k)}_{0} =
\sum_{m=0}^{\infty}
{\mathop{\sum}_{\substack{r_{1},r_{2},r_{3} \in \ZZZ_{+} \\
r_{1}+r_{2}+r_{3}=m}}}^{\hskip-.6truecm *} \hskip.6truecm
\sum_{\substack{
\n_{0},\s_{0},\n_{1},\ldots,\n_{r_{1}+r_{2}} \in\ZZZ \\
\n_{0}p+\s_{0}q+\n_{1}+\ldots+\n_{r_{1}+r_{2}}=0}}
\frac{(i\n_{0})^{r_{1}}}{r_{1}!} 
{\rm e}^{i\s_{0} t_{0}} \nonumber \\
& & \null \hskip-1.4truecm
\qquad \frac{\partial_{2}^{r_{2}} \partial_{3}^{r_{3}}}{r_{2}!r_{3}!}
G_{\n_{0}}(A_{0},C_{0}(t_{0})) \sum_{k_{1}+\ldots+k_{m}=k}
\al^{(k_{1})}_{\n_{1}}\ldots \al^{(k_{r_{1}})}_{\n_{r_{1}}} \,
A^{(k_{r_{1}+1})}_{\n_{r_{1}+1}}
\ldots A^{(k_{r_{1}+r_{2}})}_{\n_{r_{1}+r_{2}}} \,
C^{(k_{r_{1}+r_{2}+1})}\ldots C^{(k_{m})} ,
\label{eq:4.10} \end{eqnarray}
where $*$ means that the term with $r_{1}=r_{2}=0$ and $r_{3}=1$
has to be discarded --- cf. (\ref{eq:3.12}).

Therefore we see from the first equation in (\ref{eq:4.7})
that it is convenient to fix
\begin{equation}
\al_{k}(t_{0}) = \sum_{\substack{\n\in\ZZZ \\ \n\neq 0}}
\frac{F^{(k-1)}_{\n}}{i\om\nu} + \om'(A_{0})
\sum_{\substack{\n\in\ZZZ \\ \n\neq 0}}
\frac{G^{(k-1)}_{\n}}{(i\om\nu)^{2}}
\qquad \Longrightarrow \qquad \al^{(k)}_{0} = 0 ,
\label{eq:4.11} \end{equation}
so that only the functions $A^{(k)}(t)$ have the zeroth
Fourier coefficient.

In particular for $k=1$ we find
\begin{eqnarray}
& & \al^{(1)}_{\n} = \frac{1}{i\om\n}
\sum_{\substack{\n_{0},\s_{0}\in\ZZZ \\ \n_{0}p+\s_{0}q=\n}}
{\rm e}^{i\s_{0}t_{0}} F_{\n_{0},\s_{0}}(A_{0},C_{0}(t_{0})) +
\frac{\om'(A_{0})}{(i\om\n)^{2}}
\sum_{\substack{\n_{0},\s_{0}\in\ZZZ \\ \n_{0}p+\s_{0}q=\n}}
{\rm e}^{i\s_{0}t_{0}} G_{\n_{0},\s_{0}}(A_{0},C_{0}(t_{0})) ,
\nonumber \\
& & A^{(1)}_{\n} = \frac{1}{i\om\n}
\sum_{\substack{\n_{0},\s_{0}\in\ZZZ \\ \n_{0}p+\s_{0}q=\n}}
{\rm e}^{i\s_{0}t_{0}} G_{\n_{0},\s_{0}}(A_{0},C_{0}(t_{0})) ,
\label{eq:4.12} \end{eqnarray}
for $\n\neq 0$, and
\begin{equation}
%\begin{eqnarray}
%& & \al^{(1)}_{0} =  -
%\sum_{\substack{\n\in\ZZZ \\ \n\neq 0}}
%\sum_{\substack{\n_{0},\s_{0}\in\ZZZ \\ \n_{0}p+\s_{0}q=\n}}
%{\rm e}^{i\s_{0}t_{0}}
%\frac{F_{\n_{0},\s_{0}}(A_{0},C_{0}(t_{0}))}{i\om\n} -
%\sum_{\substack{\n\in\ZZZ \\ \n\neq 0}}
%\sum_{\substack{\n_{0},\s_{0}\in\ZZZ \\ \n_{0}p+\s_{0}q=\n}}
%{\rm e}^{i\s_{0}t_{0}}
%\frac{G_{\n_{0},\s_{0}}(A_{0},C_{0}(t_{0}))}{(i\om\n)^{2}} ,
%\nonumber \\
%& &
A^{(1)}_{0} = - \frac{1}{\om'(A_{0})}
\sum_{\substack{\n_{0},\s_{0}\in\ZZZ \\ \n_{0}p+\s_{0}q=0}}
{\rm e}^{i\s_{0}t_{0}}
F_{\n_{0},\s_{0}}(A_{0},C_{0}(t_{0})) ,
%\label{eq:4.12} \end{eqnarray}
\label{eq:4.13} \end{equation}
for $\n=0$, while by writing 
\begin{eqnarray}
& & C^{(1)} = - \frac{1}{D(t_{0})} \Big(
\sum_{\substack{\n_{1},\n_{2} \in\ZZZ \\ \n_{1}+\n_{2}= 0}}
\sum_{\substack{\n_{0},\s_{0}\in\ZZZ \\ \n_{0}p+\s_{0}q=\n_{1}}}
{\rm e}^{i\s_{0}t_{0}} i\n_{0} \,
G_{\n_{0},\s_{0}}(A_{0},C_{0}(t_{0})) \, \al^{(1)}_{\n_{2}}
\nonumber \\
& & \qquad \qquad \qquad +
\sum_{\substack{\n_{1},\n_{2} \in\ZZZ \\ \n_{1}+\n_{2}= 0}}
\sum_{\substack{\n_{0},\s_{0}\in\ZZZ \\ \n_{0}p+\s_{0}q=\n_{1}}}
{\rm e}^{i\s_{0}t_{0}} \partial_{2}
G_{\n_{0},\s_{0}}(A_{0},C_{0}(t_{0})) \, A^{(1)}_{\n_{2}} \Big)
\equiv C_{1}(t_{0}) ,
\label{eq:4.14} \end{eqnarray}
we can express $C^{(1)}$ in terms of the quantities in (\ref{eq:4.12}).

In order to study the convergence of the series it is convenient
to express all quantities in terms of trees. The strategy
is very simple: one iterates the relations
(\ref{eq:4.6}), (\ref{eq:4.7}) and (\ref{eq:4.9}),
which express the coefficients of order $k$ in terms
of the coefficients of lower order, until we are left only
with the coefficients of first order, for which
the explicit expressions (\ref{eq:4.12}),
(\ref{eq:4.13}) and (\ref{eq:4.14}) are at our disposal.

Trees are defined in the standard way. We briefly recall
the basic notations, by referring to \cite{GBG}
for an introductory review and further details, and also
to \cite{GG1,GBD1} for a discussion in similar contexts.

A tree $\theta$  is defined as a partially ordered set of points,
connected by oriented \textit{lines}. The lines are consistently
oriented toward a unique point $\gotr$ called the {\it root}.
The root admits only one entering line called the \textit{root line}.
All points except the root are called \textit{nodes}.
Denote with $V(\theta)$ and $L(\theta)$ the set of nodes and lines
in $\theta$, respectively, and with $|L(\theta)|$ and
$|V(\theta)|$ the number of lines and nodes of $\theta$, respectively.

If a line $\ell$ connects two points $\gotv_{1},\gotv_{2}$ and is
oriented from $\gotv_2$ to $\gotv_1$, we say that $\gotv_{2}
\prec \gotv_{1}$ and we shall write $\ell_{\gotv_{2}}=\ell$.
We shall say also that $\ell$ exits $\gotv_{2}$ and enters $\gotv_{1}$.
It can be convenient to imagine that the line $\ell$
carries an arrow pointing toward the node $\gotv_{1}$:
the arrow will be thought of as superimposed on the line itself.

More generally we write $\gotv_{2} \prec \gotv_{1}$ if $\gotv_{1}$
is on the path of lines connecting $\gotv_{2}$ to the root:
hence the orientation of the lines is opposite to the
partial ordering relation $\prec$. Along the path
from $\gotv_{2}$ to $\gotv_{1}$ all arrows point toward
$\gotv_{1}$. In particular all arrows point toward the root.

Each line $\ell$ carries a pair of labels $(h_{\ell},\de_{\ell})$,
with $h_{\ell}\in\{\al,A,C\}$ and $\de_{\ell}\in\{1,2\}$ such that
$\de_{\ell}=1$ for $h_{\ell}\neq\al$. We call $h_{\ell}$ and
$\de_{\ell}$ the \textit{component label} and
the \textit{degree label} of the line $\ell$, respectively.
Given a node $\gotv$ call $r_{\gotv 1}$, $r_{\gotv 2}$,
and $r_{\gotv 3}$ the number of lines entering $\gotv$ carrying
a component label $h=\al$, $h=A$, and $h=C$, respectively.
Hence, the values of $r_{\gotv 1},r_{\gotv 2},r_{\gotv 3}$
are uniquely determined by the component labels of
the lines entering $\gotv$.

We associate with each node $\gotv$ two \textit{mode labels}
$\n_{\gotv},\s_{\gotv}\in\ZZZ$ and we also set for convenience
$\de_{\gotv}= \de_{\ell_{\gotv}}$. With each line $\ell$ we
associate a further label $\n_{\ell} \in \ZZZ$,
called the \textit{momentum} of the line, such that
\begin{equation}
\n_{\ell} = \n_{\ell_{\gotv}} =
\sum_{\substack{\gotw \in V(\theta) \\ \gotw \preceq \gotv }}
\left( \n_{\gotw} + \s_{\gotw} \right) ,
\label{eq:4.15} \end{equation}
with the constraints that $\n_{\ell}=0$ if $h_{\ell}=C$ 
and $\n_{\ell} \neq 0$ if $h_{\ell}=\al$. The relation (\ref{eq:4.15})
expresses a conservation law at each node: the momentum
of the line exiting $\gotv$ is the sum of the
momenta of the lines entering $\gotv$ plus
the mode labels of the node $\gotv$ itself.
Note that the momentum ``flows'' through each line
in the sense of the arrow superimposed on the line.

The trees with all the labels listed above are
called \textit{labelled trees}. Then given a labelled tree $\theta$
we associate with each line $\ell$ a \textit{propagator}
\begin{equation}
g_{\ell} = \begin{cases}
{\displaystyle \frac{\om'(A_{0})^{\de_{\ell}-1}}{(i\om
\n_{\ell})^{\de_{\ell}}}} , &
h_{\ell} = \al, A , \quad \n_{\ell} \neq 0 , \\ \\
{\displaystyle - \frac{1}{\om'(A_{0})} } , &
h_{\ell} = A , \quad \n_{\ell} = 0 , \\ \\
{\displaystyle - \frac{1}{D(t_{0})} } , &
h_{\ell} = C , \quad \n_{\ell} = 0 ,
\end{cases}
\label{eq:4.16} \end{equation}
and with each node $\gotv$ a \textit{node factor}
\begin{equation}
N_{\gotv} = \begin{cases}
{\displaystyle
\frac{(i\n_{0})^{r_{\gotv 1}}\partial_{2}^{r_{\gotv 2}}
\partial_{3}^{r_{\gotv 3}}}{r_{\gotv 1}!r_{\gotv 2}!r_{\gotv 3}!}
{\rm e}^{i\s_{\gotv} t_{0}} F_{\n_{\gotv}}(A_{0},C_{0}(t_{0})) } , &
h_{\gotv} = \al , \quad \de_{\gotv}=1 , \\ \\
{\displaystyle
\frac{(i\n_{0})^{r_{\gotv 1}}\partial_{2}^{r_{\gotv 2}}
\partial_{3}^{r_{\gotv 3}}}{r_{\gotv 1}!r_{\gotv 2}!r_{\gotv 3}!}
{\rm e}^{i\s_{\gotv} t_{0}} G_{\n_{\gotv}}(A_{0},C_{0}(t_{0})) } , &
h_{\gotv} = \al , \quad \de_{\gotv}=2 , \\ \\
{\displaystyle
\frac{(i\n_{0})^{r_{\gotv 1}}\partial_{2}^{r_{\gotv 2}}
\partial_{3}^{r_{\gotv 3}}}{r_{\gotv 1}!r_{\gotv 2}!r_{\gotv 3}!}
{\rm e}^{i\s_{\gotv} t_{0}} G_{\n_{\gotv}}(A_{0},C_{0}(t_{0})) } , &
h_{\gotv} = A , \quad \de_{\gotv}=1 , \\ \\
{\displaystyle
\frac{(i\n_{0})^{r_{\gotv 1}}\partial_{2}^{r_{\gotv 2}}
\partial_{3}^{r_{\gotv 3}}}{r_{\gotv 1}!r_{\gotv 2}!r_{\gotv 3}!}
{\rm e}^{i\s_{\gotv} t_{0}} G_{\n_{\gotv}}(A_{0},C_{0}(t_{0})) } , &
h_{\gotv} = C , \quad \de_{\gotv}=1 ,
\end{cases}
\label{eq:4.17} \end{equation}
with the constraint that when $h_{\gotv}=C$ (and $\de_{\gotv}=1$)
one has either $r_{\gotv 3} \ge 2$ or $r_{\gotv 1}+r_{\gotv 2} \ge 1$.
This constraint reflects the condition $*$ in (\ref{eq:4.10}).

Finally we define the \textit{value} of a tree $\theta$ the number
\begin{equation}
\Val(\theta) = \Big( \prod_{\ell\in L(\theta)} g_{\ell} \Big)
\Big( \prod_{\gotv \in V(\theta)} N_{\gotv} \Big) ,
\label{eq:4.18} \end{equation}
which is a well-defined quantity: indeed all propagators
and node factors are bounded quantities.

Call the \textit{order} of the tree $\theta$ the number
\begin{equation}
k(\theta) = \left\{ \ell \in L(\theta) : h_{\ell} \neq C \right\} ,
\label{eq:4.19} \end{equation}
the \textit{total momentum} of $\theta$ the momentum $\n(\theta)$
of the root line, and the \textit{total component label} of $\theta$
the component label $h(\theta)$ associated to the root line.
The number of nodes (and lines) of any tree $\theta$ is
related to its order $k(\theta)$ as follows.

%%%%%%%%%%%%%%%%%%%%%%%%%%%%%%%%%%%%%%%%%%%%%%%%%%%%%%%%%%%%%%%%%%%%%%%%%
\begin{lemma} \label{lem:1}
For any tree $\theta$ one has $|L(\theta)|=|V(\theta)|\le 2k(\theta)$.
\end{lemma}
%%%%%%%%%%%%%%%%%%%%%%%%%%%%%%%%%%%%%%%%%%%%%%%%%%%%%%%%%%%%%%%%%%%%%%%%%

\prova The equality $|L(\theta)|=|V(\theta)|$ is obvious by construction.
We prove by induction on $k$ the bounds
\begin{equation}
|V(\theta)| \le \begin{cases}
3 k(\theta) - 2 , & h(\theta) = \al, A , \\
3 k(\theta) - 1 , & h(\theta) = C .
\end{cases}
\label{eq:4.20} \end{equation}
For $k=1$ the bound (\ref{eq:4.20}) is trivially satisfied, as a direct
check shows: simply compare (\ref{eq:4.12}) to (\ref{eq:4.14})
with the definition of trees in that case. Assume that the bound
holds for all $k'<k$, and let us show that then it holds also for
$k$. Call $\ell_{0}$ the root line of $\theta$ and
$\gotv_{0}$ the node which the root line exits. Call
$r_{1}$, $r_{2}$, and $r_{3}$ the number of lines entering
$\gotv_{0}$ with component labels $\al$, $A$, and $C$, respectively,
and denote with $\theta_{1},\ldots,\theta_{r_{1}+r_{2}+r_{3}}$
the subtrees which have those lines as root lines. Then
\begin{equation}
|V(\theta)| = 1 + \sum_{r=j}^{r_{1}+r_{2}+r_{3}} |V(\theta_{j})| .
\label{eq:4.21} \end{equation}
Then if $\ell_{0}$ has component label $h_{\ell_{0}} \in\{\al,A\}$
we have
\begin{equation}
|V(\theta)| \le 1 + 3 \left( k - 1 \right) - r_{3} - 2 \left(
r_{1}+r_{2} \right) \le 3k-3 < 3k-2  ,
\label{eq:4.22} \end{equation}
by the inductive hypothesis and by the fact that
$k(\theta_{1})+\ldots+k(\theta_{r_{1}+r_{2}+r_{3}})=k-1$,
whereas if $\ell_{0}$ has component label $h_{\ell_{0}}=C$ we have
\begin{equation}
|V(\theta)| \le 1 + 3 k - r_{3} - 2 \left(
r_{1}+r_{2} \right) \le 3k -1 ,
\label{eq:4.23} \end{equation}
by the inductive hypothesis, by the fact that
$k(\theta_{1})+\ldots+k(\theta_{r_{1}+r_{2}+r_{3}})=k$,
and by the constraint that either $r_{3}\ge2$ or $r_{1}+r_{2}\ge 1$
in such a case --- cf. the comment after (\ref{eq:4.17}).
Therefore the assertion is proved. \qed

\*

Define $\Theta_{k,\n,h}$ the set of all trees
of order $k(\theta)=k$, total momentum $\n(\theta)=\n$,
and total component label $h(\theta)=h$.
By collecting together all the definitions given above,
one obtains that the Fourier coefficients $\al^{(k)}_{\n}$
and $A^{(k)}_{\n}$ and the constants $C_{k}$ can be written
for all $k\ge 1$ in terms of trees as
\begin{eqnarray}
& & \al^{(k)}_{\n} = \sum_{\theta \in \Theta_{k,\n,\al}}
\Val(\theta) , \qquad \n \neq 0 , \qquad \qquad
\al^{(k)}_{0} = 0 \nonumber \\
& & A^{(k)}_{\n} = \sum_{\theta \in \Theta_{k,\n,A}}
\Val(\theta) , \qquad
C^{(k)} = \sum_{\theta \in \Theta_{k,0,C}} \Val(\theta) .
\label{eq:4.24} \end{eqnarray}
The proof of (\ref{eq:4.24}) can be performed by induction; cf.
\cite{GBG} for details.

The number of unlabelled trees of order $k$ is bounded
by the number of random walks of $2k$ steps, hence by $2^{2k}$ \cite{GM}.
The sum over all labels except the mode labels and the momenta
is bounded again by a constant to the power $k$ --- simply
because all such labels can assume only a finite number of values.
Finally the sum over the mode labels --- which uniquely determine
the momenta through the relation (\ref{eq:4.15}) ---
can be performed by using for each node half the exponential decay
factor ${\rm e}^{-\ka(|\n_{\gotv}|+|\s_{\gotv}|)}$
provided by the bounds (\ref{eq:4.4}).
The conclusion is that we obtain eventually the bounds
\begin{equation}
\left| \al^{(k)}_{\n} \right| \le B_{1} B_{2}^{k}
{\rm e}^{-\ka|\n|/2} , \qquad
\left| A^{(k)}_{\n} \right| \le B_{1} B_{2}^{k}
{\rm e}^{-\ka|\n|/2} ,\qquad
\left| C^{(k)} \right| \le B_{1} B_{2}^{k} , \qquad
\label{eq:4.25} \end{equation}
for suitable constants $B_{1}$ and $B_{2}$.
This proves the convergence of the series (\ref{eq:3.1})
and (\ref{eq:3.2}) for $|\eps|<\eps_{0}$,
with $\eps_{0}$ small enough.
Note that with respect to \cite{GBG} here
the analysis is much easier as there is no small divisors problem.

The construction described above provides also a useful algorithm
which can be implemented numerically in order
to compute the solution to any prescribed accuracy
(provided $\eps$ is small enough).

Now, we come back to the problem of determining the boundary
of the set in the plane $(\eps,\g)$, with $\g=\eps C$,
in which there are subharmonic solutions of order $q/p$.

We have to find the solutions of (\ref{eq:4.1}), that is,
solve the equation
\begin{equation}
0 = \frac{\partial}{\partial t_{0}} C(\eps,t_{0})
= C_{0}'(t_{0}) + \eps C_{1}'(t_{0}) +
\eps^{2} C_{2}'(t_{0}) + \ldots ,
\label{eq:4.26} \end{equation}
where $C_{k}'(t_{0})={\rm d} C_{k}(t_{0})/{\rm d}t_{0}$.

The function $t_{0} \to C(\eps,t_{0})$ is analytic
in $t_{0}$ for all $|\eps|<\eps_{0}$ (for which it is defined
and analytic in $\eps$), so that for fixed $\eps$ the
equation (\ref{eq:4.26}) can always be solved. It has at least
the two solutions $t_{0}=\tau_{1}(\eps)$ and $t_{0}=\tau_{2}(\eps)$
corresponding to the absolute minimum and to the absolute maximum,
respectively, of the function $C(\eps,t_{0})$. In general these
solutions are not smooth in $\eps$. This proves Theorem \ref{thm:1}.

Suppose now that at the value $t_{0}$ such that $C_{0}'(t_{0})=0$
one has furthermore $C_{0}''(t_{0}) \neq 0$. In that case,
if $\tau_{0}=\tau_{0}(\eps)$ is a solution of (\ref{eq:4.26}) ---
$\tau_{0}$ is a point of minimum or maximum for $C(\eps,t_{0})$ ---
then $\tau_{0}$ must be analytically close to $t_{0}$. Hence  
$\eps \to \tau_{0}(\eps)$ is an analytic function of $\eps$,
so that also $\eps \to C_{1}(\eps)$ and $\eps\to C_{2}(\eps)$
are smooth (in fact analytic) in $\eps$.
Therefore also Theorem \ref{thm:2} follows.

The last observation suggests how to extend Theorem \ref{thm:2}
to obtain smooth bifurcation curves when Hypothesis \ref{hyp:3}
fails to be satisfied.

%%%%%%%%%%%%%%%%%%%%%%%%%%%%%%%%%%%%%%%%%%%%%%%%%%%%%%%%%%%%%%%%%%%%%%%%%
\begin{hypothesis} \label{hyp:4}
There exists $k \ge 1$ such that the functions $C_{p}(t_{0})$
are identically constant in $t_{0}$ for all $p=0,\ldots,k-1$.
If $t_{m}$ and $t_{M}$ are the values in $[0,2\pi)$
for which the function $C_{k}(t_{0})$ attains its minimum and
its maximum, respectively, then $C_{k}''(t_{m}) C_{k}''(t_{M}) \neq 0$.
\end{hypothesis}
%%%%%%%%%%%%%%%%%%%%%%%%%%%%%%%%%%%%%%%%%%%%%%%%%%%%%%%%%%%%%%%%%%%%%%%%%

The following result extends Theorem \ref{thm:2},
as it deals with the case in which the subharmonic Melnikov function
does not depend explicitly on $t_{0}$, that is $C_{0}'(t_{0})\equiv0$.

%%%%%%%%%%%%%%%%%%%%%%%%%%%%%%%%%%%%%%%%%%%%%%%%%%%%%%%%%%%%%%%%%%%%%%%%%
\begin{theorem} \label{thm:3}
Consider the system (\ref{eq:2.1}) and assume that Hypotheses
\ref{hyp:1}, \ref{hyp:2} and \ref{hyp:4} hold for the resonant torus
with energy $A_{0}$ such that $\om(A_{0})=p/q$. There exist
$\eps_{0} > 0$ and two functions $\widetilde \g_{1}(\eps)$
and $\widetilde \g_{2}(\eps)$, analytic for $|\eps|<\eps_{0}$,
with $\widetilde\g_{1}(0)=
\widetilde\g_{2}(0)$ and $\widetilde\g_{1}(\eps)\neq \widetilde
\g_{2}(0)$ for all $\eps\neq0$, such that the two functions
\begin{equation}
\g_{1}(\eps) = \begin{cases}
\max\{\widetilde \g_{1}(\eps),\widetilde \g_{2}(\eps) \} ,
& \eps > 0 , \\
\min\{\widetilde \g_{1}(\eps),\widetilde \g_{2}(\eps) \} ,
& \eps < 0 , \end{cases}
\qquad\qquad
\g_{2}(\eps) = \begin{cases}
\min\{\widetilde \g_{1}(\eps),\widetilde \g_{2}(\eps) \} ,
& \eps > 0 , \\
\max\{\widetilde \g_{1}(\eps),\widetilde \g_{2}(\eps) \} ,
& \eps < 0 , \end{cases}
\label{eq:4.27} \end{equation}
have the same tangent lines at the origin, and (\ref{eq:2.1})
has at least one subharmonic solution of order $q/p$ for $\g_{2}(\eps)
\le \eps C \le \g_{1}(\eps)$ when $\eps\in(0,\eps_{0})$ and for
$\g_{1}(\eps)\le \eps C \le \g_{2}(\eps)$ when $\eps\in(-\eps_{0},0)$.
\end{theorem}
%%%%%%%%%%%%%%%%%%%%%%%%%%%%%%%%%%%%%%%%%%%%%%%%%%%%%%%%%%%%%%%%%%%%%%%%%

The proof follows the same lines as that of Theorem \ref{thm:2}.
The only difference is that up to order $k-1$ the initial
phase is left undetermined. In fact to first order one has
$M(t_{0},C)=M(C)=0$ which fixes $C=C_{0}$ (by Hypothesis \ref{hyp:2}),
while to orders $k'=2,\ldots,k-1$ the constants $C_{k}$
are fixed and are independent of $t_{0}$ by Hypothesis \ref{hyp:4}.
Then we can write $C(\eps,t_{0})=\gotC_{1}(\eps)+\gotC_{2}(\eps,t_{0})$,
with $\gotC_{1}(\eps)=C_{0}+\eps C_{1} + \ldots + \eps^{k-1} C_{k-1}$
and $\gotC_{2}(\eps,t_{0})=\eps^{k}(C_{k}(t_{0})+O(\eps))$, and
from order $k$ on the constants $C_{k}$ are fixed as functions
of $t_{0}$. Moreover equation (\ref{eq:4.26}) reduces to
$0=C_{k}'(t_{0}) + \eps C_{k+1}'(t_{0})+\ldots$. Therefore
we can reason as in the previous case ($k=0$) and
we find that $C_{k}(t_{0})$ has at least two stationary
points $t_{0}=t_{1}$ and $t_{0}=t_{2}$, corresponding
to the minimum point and to the maximum point, respectively.
By Hypothesis \ref{hyp:4} also $\gotC_{2}(\eps,t_{0})$ has two
stationary points at $\tau_{1}(\eps)=t_{1}+O(\eps)$ and
$\tau_{2}(\eps)=t_{2}+O(\eps)$, with $\tau_{1}(\eps)$ and
$\tau_{2}(\eps)$ analytic in $\eps$ for $\eps$ small enough.
Then we can define $\widetilde \g_{1}(\eps)=C(\eps,\tau_{1}(\eps))$
and $\widetilde \g_{2}(\eps)=C(\eps,\tau_{2}(\eps))$:
by construction, both $\widetilde \g_{1}(\eps)$ and $\widetilde
\g_{1}(\eps)$ are analytic in $\eps$ for $\eps$ small enough.
If we define $\g_{1}(\eps)$ and $\g_{2}(\eps)$ according
to (\ref{eq:4.27}) then the proof of the theorem is achieved.

Note that in this case the definition (\ref{eq:4.27}) coincides
with the general definition (\ref{eq:4.1}) for the bifurcation curves.
Furthermore, if we assume Hypothesis \ref{hyp:3} instead of
Hypothesis \ref{hyp:4}, then one has $\widetilde \g_{1}(\eps)=
\g_{1}(\eps)$ and $\widetilde \g_{1}(\eps)=\g_{1}(\eps)$, so that
also $\g_{1}(\eps)$ and $\g_{2}(\eps)$ are analytic,
as stated in Theorem \ref{thm:2}.

Finally we note that if the functions $C_{k}(t_{0})$ are identically
constant in $t_{0}$ for all $k\in\ZZZ_{+}$ then one has $C(\eps,t_{0})=
C(\eps)$. In this case the two curves $\g_{1}(\eps)$ and
$\g_{2}(\eps)$ coincide, and all values of $t_{0}$ are allowed.
This means that the whole manifold corresponding to the resonant torus
persists. On the other hand the parameter $C$ must be fixed in a very
precise way, as a function of $\eps$, and any small deviation from
that value destroys the torus. This result can be compared
with \cite{C,CL}, where a similar situation is discussed.

For $(\eps,\g)$ inside the set of existence of subharmonic solutions
one can investigate how many of them exist. For $p=1$ the initial
phase $t_{0}$ varies in the interval $[0,2\pi q]$, where $T_{0}=2\pi q$
is the period of the unperturbed periodic solution. The function
$C(\eps,t_{0})$ has period $2\pi$ in $t_{0}$, so that
it is repeated $q$ times in the interval $[0,2\pi q]$.
Hence for any fixed value $|\eps|<\eps_{0}$ and any $C$ strictly
between the maximum and the mininum value attained by the function
$t_{0} \to C(\eps,t_{0})$ there are at least $2q$ values $t_{i}$,
$i=2,\ldots,2q$, such that $C=C(\eps,t_{i})$. If $C$ coincides with
either its maximum or its minimum then there are at least $q$ values
$t_{i}$, $i=2,\ldots,q$, such that $C=C(\eps,t_{i})$. Therefore we can
conclude that, for $p=1$, inside the set of existence of subharmonic
solutions there are at least $2q$ such solutions, as found
in \cite{CH}, while on the boundary of that set there are $q$ of them.

We can summarise the discussion above in the following statement.

%%%%%%%%%%%%%%%%%%%%%%%%%%%%%%%%%%%%%%%%%%%%%%%%%%%%%%%%%%%%%%%%%%%%%%%%%
\begin{theorem} \label{thm:4}
Under the same assumptions of Theorem \ref{thm:1} assume $p=1$.
Take $|\eps|<\eps_{0}$, and for such values of $\eps$ let
$\eps\to\g_{1}(\eps)$ and $\eps\to\g_{2}(\eps)$ be the two
bifurcation curves whose existence is assured by Theorem \ref{thm:1}.
For $\min\{\g_{1}(\eps),\g_{2}(\eps)\} < \g < \max\{\g_{1}(\eps),
\g_{2}(\eps)\}$ there at least $2q$ subharmonic solutions of order $q$.
If either $\g=\g_{1}(\eps)$ or $\g=\g_{2}(\eps)$ one has
at least $q$ subharmonic solutions of order $q$.
\end{theorem}
%%%%%%%%%%%%%%%%%%%%%%%%%%%%%%%%%%%%%%%%%%%%%%%%%%%%%%%%%%%%%%%%%%%%%%%%%

Theorem \ref{thm:4} should be compared with Theorem 2.1 in \cite{CH}.

%%%%%%%%%%%%%%%%%%%%%%%%%%%%%%%%%%%%%%%%%%%%%%%%%%%%%%%%%%%%%%%%%%%%%%%%%
%%%%%%%%%%%%%%%%%%%%%%%%%%%%%%%%%%%%%%%%%%%%%%%%%%%%%%%%%%%%%%%%%%%%%%%%%
\zerarcounters
\section{Application to dissipative systems with forcing}
\label{sec:5}
%%%%%%%%%%%%%%%%%%%%%%%%%%%%%%%%%%%%%%%%%%%%%%%%%%%%%%%%%%%%%%%%%%%%%%%%%
%%%%%%%%%%%%%%%%%%%%%%%%%%%%%%%%%%%%%%%%%%%%%%%%%%%%%%%%%%%%%%%%%%%%%%%%%

Let us consider a one-dimensional system, subject
to a conservative force $g(x)$, in the presence of dissipation
and of a periodic forcing. If the periodic forcing and the
dissipation coefficient are both small we can write the equations 
for the system as
\begin{equation}
\ddot x + g(x) +\g \dot x = \eps f(x,t) ,
\qquad \g = \eps \, C ,
\label{eq:5.1} \end{equation}
where $\eps f(x,t)$ is the forcing of period $2\pi$ and
$C$ is a parameter. Assume that both $g$ and $f$ are analytic
in their arguments. If $f$ depends only on $t$, equation
(\ref{eq:5.1}) reduces to the equation studied in \cite{HT1}.

Let us assume that the unperturbed system ($\eps=0$)
is Liouville-integrable and anysochronous. This means that,
in action-angle variables, the equations (\ref{eq:5.1})
can be written in the form (\ref{eq:2.1}), and,
furthermore, that Hypothesis \ref{hyp:1} is satisfied.

We define the subharmonic Melnikov function in terms of the action-angle
variable as in (\ref{eq:2.2}). To check that Hypothesis
\ref{hyp:2} is also satisfied we use the following result.

%%%%%%%%%%%%%%%%%%%%%%%%%%%%%%%%%%%%%%%%%%%%%%%%%%%%%%%%%%%%%%%%%%%%%%%%%
\begin{lemma} \label{lem:2}
The subharmonic Melnikov function is invariant under a
transformation of coordinates.
\end{lemma}
%%%%%%%%%%%%%%%%%%%%%%%%%%%%%%%%%%%%%%%%%%%%%%%%%%%%%%%%%%%%%%%%%%%%%%%%%

\prova Consider a system of differential equations in $\RRR^{2}$
\begin{equation}
\dot x = f(x) + \eps g(x,t) ,
\label{eq:5.2} \end{equation}
and define the subharmonic Melnikov function \cite{M,GH,CH} for
a subharmonic solution $x_{0}(t)$ of period $T$ as
\begin{equation}
M(t_{0}) = \frac{1}{T} \int_{0}^{T} 
{\rm d}t \big( f_{1}(x_{0}(t)) \, g_{2}(x_{0}(t)) -
f_{2}(x_{0}(t)) \, g_{1}(x_{0}(t)) \big).
\label{eq:5.3} \end{equation}
Take the transformation of coordinates $\xi \to x = h(\xi)$.
In the new coordinates the system reads
\begin{equation}
\dot \xi = \phi(\xi) + \eps \g(\xi,t) ,
\label{eq:5.4} \end{equation}
where $\phi(\xi)=\partial h^{-1}(h(\xi)) \, f(h(\xi))$ and
$\g(\xi)=\partial h^{-1}(h(\xi)) \, g(h(\xi))$,
and the subharmonic Melnikov function becomes
\begin{equation}
\MM(t_{0}) = \frac{1}{T} \int_{0}^{T}  {\rm d}t
\big( \phi_{1}(\xi_{0}(t)) \, \g_{2}(\xi_{0}(t)) -
\phi_{2}(\xi_{0}(t)) \, \g_{1}(\xi_{0}(t)) \big) ,
\label{eq:5.5} \end{equation}
where $\xi_{0}(t)$ is the subharmonic solution expressed
in the new variables.

By noting that
\begin{equation}
\partial h^{-1} (h(\xi)) =
\left( \partial h(\xi) \right)^{-1} =
\frac{1}{J} \left( \begin{matrix}
\partial_{2} h_{2}(\xi) & - \partial_{2} h_{1}(\xi) \\
- \partial_{1} h_{2}(\xi) & \partial_{1} h_{1}(\xi)
\end{matrix} \right) ,
\label{eq:5.6} \end{equation}
where $J=\det \partial h = \partial_{1} h_{1} \partial_{2} h_{2} -
\partial_{1} h_{2} \partial_{2} h_{1}$ is the Jacobian
of the transformation, one obtains
\begin{eqnarray}
& & \MM(t_{0}) = \frac{1}{T} \int_{0}^{T} {\rm d}t \frac{1}{J}
\Big(
\left( \partial_{2} h_{2} f_{1} - \partial_{2} h_{1} f_{2} \right)
\left( -\partial_{1} h_{2} g_{1} + \partial_{1} h_{1} g_{2} \right) -
\nonumber \\
& & \qquad \qquad \qquad
\left( -\partial_{1} h_{2} f_{1} + \partial_{1} h_{1} f_{2} \right)
\left( \partial_{2} h_{2} g_{1} - \partial_{2} h_{1} g_{2} \right)
\Big) \nonumber \\
& & \qquad \quad \,
= \frac{1}{T} \int_{0}^{T} {\rm d}t \frac{1}{J}
\left( \partial_{1} h_{1} \partial_{2} h_{2} -
\partial_{1} h_{2} \partial_{2} h_{1} \right)
\left( f_{1} g_{2} - f_{2} g_{1} \right) ,
\label{eq:5.7} \end{eqnarray}
where the function $h$ is computed in $\xi_{0}(t)$
and the functions $f,g$ are computed in $x_{0}(t)=h(\xi_{0}(t))$.
Hence (\ref{eq:5.3}) yields $\MM(t_{0})=M(t_{0})$,
so that the assertion follows. \qed

\*

This means that we can compute the subharmonic Melnikov function
for the system (\ref{eq:5.1}) in the coordinates
$(x,y) = (x,\dot x)$.
In that case the unperturbed vector field is
$(y,-g(x))$ and the perturbation reads $(0,-\eps C y + \eps f(x,t))$,
so that the subharmonic Melnikov function becomes
\begin{equation}
M(t_{0},C) = \frac{1}{T} \int_{0}^{T} {\rm d}t \,
y_{0}(t) \big( - C y_{0}(t) + f(x_{0}(t),t+t_{0}) \big) =
- C \langle y_{0}^{2} \rangle +
\langle y_{0} f(x_{0}(\cdot),\cdot + t_{0}) \rangle .
\label{eq:5.8} \end{equation}
Therefore the subharmonic Melnikov function vanishes provided
$C=C_{0}(t_{0})$, where $C_{0}(t_{0}) = (\langle y_{0}^{2} \rangle)^{-1}
\langle y_{0} f(x_{0}(\cdot),\cdot + t_{0}) \rangle$,
which is well-defined because $\langle y_{0}^{2} \rangle >0$.
Moreover one has $\partial M(t_{0},C)/\partial C=
-\langle y_{0}^{2} \rangle \neq 0$. Therefore
Hypothesis \ref{hyp:2} is also satisfied, and
Theorem \ref{thm:2} applies to the system (\ref{eq:5.1}).

We can state our result as follows.

%%%%%%%%%%%%%%%%%%%%%%%%%%%%%%%%%%%%%%%%%%%%%%%%%%%%%%%%%%%%%%%%%%%%%%%%%
\begin{theorem} \label{thm:5}
Consider the system (\ref{eq:5.1}) and assume that Hypothesis
\ref{hyp:1} holds for the invariant torus with energy $A_{0}$
such that $\om(A_{0})=p/q$. There exist $\eps_{0} > 0$ and two
continuous functions $\g_{1}(\eps)$ and $\g_{2}(\eps)$, with
$\g_{1}(0)=\g_{2}(0)$, $\g_{1}(\eps)\ge \g_{2}(\eps)$ for $\eps\ge 0$
and $\g_{1}(\eps)\le \g_{2}(\eps)$ for $\eps\le 0$,
such that (\ref{eq:2.1}) has at least one
subharmonic solution of period $2\pi p$ for $\g_{2}(\eps) \le \eps C
\le \g_{1}(\eps)$ when $\eps\in(0,\eps_{0})$ and for
$\g_{1}(\eps)\le \eps C \le \g_{2}(\eps)$ when $\eps\in(-\eps_{0},0)$.
\end{theorem}
%%%%%%%%%%%%%%%%%%%%%%%%%%%%%%%%%%%%%%%%%%%%%%%%%%%%%%%%%%%%%%%%%%%%%%%%%

Of course Theorem \ref{thm:5} is a corollary of Theorem \ref{thm:1}.
It should be compared with Corollary 2.3 in \cite{CH}
(cf. also \cite{HT1}). Our result is stronger as it requires,
in Chow and Hale's notations, only Hypothesis (H$_{1}$),
which corresponds to our Hypothesis \ref{hyp:1}.
If one assumes also Hypothesis (H$_{4}$) of \cite{CH}, which
corresponds to our hypothesis \ref{hyp:3}, then Theorem \ref{thm:2}
applies, and the result of \cite{CH} is recovered.

One expects that, in the case of system (\ref{eq:5.1}),
the two bifurcation curves $\g_{1}(\eps)$ and $\g_{2}(\eps)$
contain the real axis, that is $\min\{\g_{1}(\eps),\g_{2}(\eps)\}
\le 0 \le \max\{\g_{1}(\eps),\g_{2}(\eps)\}$.
Indeed for $\g=0$ the equation (\ref{eq:5.1}) describes
a quasi-integrable Hamiltonian system, and existence of
periodic solutions is well known in this case, at least
under some non-degeneracy condition on the unperturbed system,
such as Hypothesis \ref{hyp:1}. If $C_{0}(t_{0})$ is not zero
then it is easy to check that the set of existence of
subharmonic solutions includes the real axis. Indeed, this
follows from the following result.

%%%%%%%%%%%%%%%%%%%%%%%%%%%%%%%%%%%%%%%%%%%%%%%%%%%%%%%%%%%%%%%%%%%%%%%%%
\begin{lemma} \label{lem:3}
The function $C_{0}(t_{0})$ has zero mean.
\end{lemma}
%%%%%%%%%%%%%%%%%%%%%%%%%%%%%%%%%%%%%%%%%%%%%%%%%%%%%%%%%%%%%%%%%%%%%%%%%

\prova Call
\begin{equation}
F(x_{0}(t)) = \int_{0}^{2\pi} \frac{{\rm d}t_{0}}{2\pi} \,
f(x_{0}(t),t+t_{0}) = \int_{0}^{2\pi} \frac{{\rm d}t_{0}}{2\pi} \,
f(x_{0}(t),t_{0}) .
\label{eq:5.9} \end{equation}
By (\ref{eq:5.8}) the mean (with respect to $t_{0}$)
of $C_{0}(t_{0})$ is
\begin{eqnarray}
& & \int_{0}^{2\pi} \frac{{\rm d}t_{0}}{2\pi} \, C_{0}(t_{0}) =
\frac{1}{\langle y_{0}^{2} \rangle}
\int_{0}^{2\pi} \frac{{\rm d}t_{0}}{2\pi}
\int_{0}^{T} \frac{{\rm d}t}{T} \, \dot x_{0}(t) f(x_{0}(t),t+t_{0})
\nonumber \\
& & \qquad \qquad \qquad \; \; \, =
\int_{0}^{T} \frac{{\rm d}t}{T} \, \dot x_{0}(t) F(x_{0}(t)) ,
\label{eq:5.10} \end{eqnarray}
which vanishes, as the integrand can be written as a total derivative
with respect to $t$. \qed

\*

In particular Lemma \ref{lem:3} implies that if $C_{0}(t_{0})$
is not identically constant then its maximum is strictly positive
and its minimum is strictly negative, hence $\max\{\g_{1}(\eps),
\g_{2}(\eps)\}>0$ and $\min\{\g_{1}(\eps),\g_{2}(\eps)\}<0$.

To extend the same result to the case in which the functions
$C_{k'}(t_{0})$ are identically constant in $t_{0}$ for all $k'\le k-1$,
with $k\ge 1$ arbitrarily high, is more delicate, and it requires
some work. One can reason as follows.

%%%%%%%%%%%%%%%%%%%%%%%%%%%%%%%%%%%%%%%%%%%%%%%%%%%%%%%%%%%%%%%%%%%%%%%%%
\begin{lemma} \label{lem:4}
Assume that for some $\Bar k\in\ZZZ_{+}$ the coefficients
$C_{k'}(t_{0})$ vanish identically for all $k'=0,\ldots,\Bar k-1$.
Then $C_{\Bar k}(t_{0})$ has zero mean in $t_{0}$.
\end{lemma}
%%%%%%%%%%%%%%%%%%%%%%%%%%%%%%%%%%%%%%%%%%%%%%%%%%%%%%%%%%%%%%%%%%%%%%%%%

\prova Write the system (\ref{eq:5.1}) in action-angle variables.
Then there exists a Hamiltonian function $H(\al,A,t,\eps)=
H_{0}(A) + \eps H_{1}(\al,A,t)$ such that
$\om(A)=\partial_{A}H_{0}(A)$ and
\begin{equation}
\begin{cases}
\dot \al = \om(A) + \eps \partial_{A} H_{1}(\al,A,C,t) +
\eps C \, \Phi(\al,A) , \\
\dot A = - \eps \partial_{\al} H_{1}(\al,A,C,t) + \eps C \, \Psi(\al,A) ,
\end{cases}
\label{eq:5.11} \end{equation}
where $\Phi=-y\,\partial \al/\partial y$ and
$\Psi=y\,\partial A/\partial y$. Then (\ref{eq:4.6}) become
\begin{eqnarray}
& & \al^{(k)}_{\n} = \frac{1}{i\om\nu}
\left( \partial_{A} H_{1}^{(k-1)} \right)_{\n} +
\om'(A_{0}) \frac{1}{(i\om\nu)^{2}}
\left( - \partial_{\al} H_{1}^{(k-1)} \right)_{\n} +
\left( C \Phi \right)^{(k-1)}_{\n} , \nonumber \\
& & A^{(k)}_{\n} = \frac{1}{i\om\n}
\left( - \partial_{\al} H_{1}^{(k-1)} \right)_{\n} +
\left( C \Psi \right)^{(k-1)}_{\n} ,
\label{eq:5.12} \end{eqnarray}
for all $k\in\NNN$ and all $\n\neq0$. Moreover (\ref{eq:4.9}) reads
\begin{equation}
\sum_{k'=0}^{k} C_{k'} \Psi^{(k')}_{0} +
\Bar \Gamma^{(k)}_{0} = 0 , \qquad \Bar \Gamma^{(k)}_{0} =
\left( - \partial_{\al} H_{1}^{(k-1)} \right)_{0} ,
\label{eq:5.13} \end{equation}
which, for $k=\Bar k$, gives $\Gamma^{(\Bar k)}_{0}=
\Bar\Gamma^{(\Bar k)}_{0}$ and $C_{\Bar k}\Psi^{(0)}_{0}+
\Bar\Gamma^{(\Bar k)}_{0}=0$ because $C_{1}=\ldots=C_{\Bar k-1}=0$
by assumption. Moreover $\Psi^{(0)}=-\langle y_{0}^{2} \rangle \neq 0$,
by Lemma \ref{lem:2} and Hypothesis \ref{hyp:2}.

Therefore $C_{k}=C^{(k)}$, with $C^{(k)}$ given by the sum
(\ref{eq:4.24}) of tree values. We can split the set $\Theta_{k,0,C}$
into the union of disjoint families $\FF$ as follows.
Given a tree $\theta \in \Theta_{k,0,C}$ call $\gotv_{0}$ the node
which is connected to the root through the root line, and define
$V_{0}(\theta)$ as the subset of nodes $\gotv\in V(\theta)$
such that all the lines $\ell$ along the path
connecting $\gotv$ to $\gotv_{0}$ have $\n_{\ell}\neq 0$.
Then define $\FF=\FF(\theta)$ as the set of trees obtained from $\theta$
by ``shifting'' the root line to any node in $V_{0}(\theta)$,
i.e. by attaching the root line to any node $\gotv\in V_{0}(\theta)$.
Of course, as a consequence of the shift of the root line
from $\gotv_{0}$ to $\gotv$, the arrows of all lines
along the path between the two nodes are reversed. If one recalls
the diagrammatic rules introduced in Section \ref{sec:4} to
associate with any tree $\theta$ a value $\Val(\theta)$, this means
that all lines with labels $(h,\delta)=(\al,1)$ are transformed
into lines with labels $(h,\delta)=(A,1)$. Moreover the momenta
of all such lines change sign. The latter property can be seen
as follows. The momentum is defined as the sum of all mode labels
of the nodes preceding the lines --- cf. (\ref{eq:4.15}) ---
and the sum of all the mode labels is zero for any tree $\theta\in
\Theta_{k,0,C}$: then, when the arrow of a line $\ell$ is reversed
the nodes preceding $\ell$ become the nodes following $\ell$
and vice versa, so that $\n_{\ell}$ becomes $-\n_{\ell}$.
Hence the propagators of the lines $\ell$ with $\delta_{\ell}=1$
change sign, whereas the propagators of the lines $\ell$ with
$\delta_{\ell}=2$ are left unchanged. As a consequence, for each tree
$\theta'\in\FF(\theta)$ we can write $\Val(\theta)= i\n_{\gotv}
\overline{\Val}(\theta)$, where $\gotv$ is the node $\gotv\in V_{0}(\theta)$
which the root line exits and $\overline{\Val}(\theta)$ is the same
quantity for all $\theta'\in\FF(\theta)$. Therefore
\begin{equation}
\sum_{\theta'\in\FF(\theta)} \Val(\theta) =
\overline{\Val}(\theta) \sum_{\gotv\in V_{0}(\theta)} i\n_{\gotv} .
\label{eq:5.14} \end{equation}
Moreover one has
\begin{equation}
\sum_{\gotv\in V(\theta)}
\left( \n_{\gotv} + \s_{\gotv} \right) = 0
\quad \Longrightarrow \quad
\sum_{\gotv\in V_{0}(\theta)}
\left( \n_{\gotv} + \s_{\gotv} \right) = 0
\quad \Longrightarrow \quad
\sum_{\gotv\in V_{0}(\theta)} \n_{\gotv} = -
\sum_{\gotv\in V_{0}(\theta)} \s_{\gotv} ,
\label{eq:5.15} \end{equation}
so that the mean in $t_{0}$ of (\ref{eq:5.14}) gives
\begin{equation}
\int_{0}^{2\pi} \frac{{\rm d}t_{0}}{2\pi} \,
\sum_{\theta'\in\FF(\theta)} \Val(\theta) =
\int_{0}^{2\pi} \frac{{\rm d}t_{0}}{2\pi} \, \overline{\Val}(\theta) 
\sum_{\gotv\in V_{0}(\theta)} i\n_{\gotv}
= - \int_{0}^{2\pi} \frac{{\rm d}t_{0}}{2\pi} \, \overline{\Val}(\theta)
\sum_{\gotv\in V_{0}(\theta)} i\s_{\gotv}  = 0 ,
\label{eq:5.16} \end{equation}
because the mean is the sum over all labels $\s_{\gotv} \in V(\theta)$
such that $\sum_{v\in V(\theta)}\s_{\gotv}=\sum_{v\in V_{0}(\theta)}
\s_{\gotv}=0$. By using the fact that the set $\Theta_{k,0,C}$ can be written as a
disjoint union of the sets $\FF$, we obtain that $\Bar\Gamma^{(\Bar
k)}_{0}$ has zero mean in $t_{0}$, so that the assertion follows. \qed

\*

%%%%%%%%%%%%%%%%%%%%%%%%%%%%%%%%%%%%%%%%%%%%%%%%%%%%%%%%%%%%%%%%%%%%%%%%%
\begin{lemma} \label{lem:5}
Assume that for some $\Bar k\in\ZZZ_{+}$ the coefficients
$C_{k'}(t_{0})$ are identically constant for all $k'=0,\ldots,\Bar k-1$.
Then $C_{k'}(t_{0})\equiv 0$ for all $k'=0,\ldots,\Bar k-1$.
\end{lemma}
%%%%%%%%%%%%%%%%%%%%%%%%%%%%%%%%%%%%%%%%%%%%%%%%%%%%%%%%%%%%%%%%%%%%%%%%%

\prova The proof is by induction. Fix $0\le k< \Bar{k}$,
and assume that $C_{k'}(t_{0})\equiv 0$ for all $k'\le k-1$.
Then by Lemma \ref{lem:4} the function $C_{k}(t_{0})$
has zero mean. Since it is constant by hypothesis then
$C_{k}(t_{0})\equiv 0$. \qed

\*

Let $k\in\ZZZ_{+}$ be  such that $C_{k'}(t_{0})$ is identically
constant in $t_{0}$ for $k'=0,\ldots,k-1$ whereas
$C_{k}(t_{0})$ depends explicitly on $t_{0}$.
If $k=0$ this simply means that $C_{0}(t_{0})$ depends explicitly
on $t_{0}$. By Lemma \ref{lem:5} one has $C_{k'}(t_{0})\equiv 0$
for all $k'\le k-1$, and by Lemma \ref{lem:4} the function
$C_{k}(t_{0})$ has zero mean in $t_{0}$. Since $C_{k}(t_{0})$
is not identically constant then $\sup_{t_{0}\in[0,2\pi)} C_{k}(t_{0})>0$
and $\inf_{t_{0}\in[0,2\pi)} C_{k}(t_{0})<0$. Furthermore, in such a case
$C(\eps,t_{0})=\eps^{k} (C_{k}(t_{0})+O(\eps))$, so that also
\begin{equation}
\sup_{t_{0}\in[0,2\pi)} C(\eps,t_{0})>0 , \qquad
\inf_{t_{0}\in[0,2\pi)} C(\eps,t_{0})<0 ,
\label{eq:5.17} \end{equation}
for $\eps$ small enough. If we recall the definition (\ref{eq:4.1})
of the bifurcation curves we can formulate the following result.

%%%%%%%%%%%%%%%%%%%%%%%%%%%%%%%%%%%%%%%%%%%%%%%%%%%%%%%%%%%%%%%%%%%%%%%%%
\begin{theorem} \label{thm:6}
Under the same assumptions of Theorem \ref{thm:5} let
$\eps\to\g_{1}(\eps)$ and $\eps\to\g_{2}(\eps)$ be the two
bifurcation curves whose existence is assured by Theorem \ref{thm:5}.
One has $\g_{1}(\eps) \ge 0 \ge \g_{2}(\eps)$ for
$\eps\in(0,\eps_{0})$ and $\g_{1}(\eps) \le 0 \le \g_{2}(\eps)$
for $\eps\in(0,\eps_{0})$.
\end{theorem}
%%%%%%%%%%%%%%%%%%%%%%%%%%%%%%%%%%%%%%%%%%%%%%%%%%%%%%%%%%%%%%%%%%%%%%%%%

As (\ref{eq:5.17}) shows, if there is $k\ge0$ such that
$C_{k'}(t_{0})\equiv 0$ for $k'=0,\ldots,k-1$ and
$C_{k}(t_{0}) \neq 0$, then one has the strict inequalities
$\g_{1}(\eps) > 0 > \g_{2}(\eps)$ for
$\eps\in(0,\eps_{0})$ and $\g_{1}(\eps) < 0 < \g_{2}(\eps)$
for $\eps\in(0,\eps_{0})$. On the contrary if all
$C_{k}$ vanish identically, so that the full function
$C(\eps,t_{0})$ has to be zero, then $\g_{1}(\eps)=\g_{2}(\eps)=0$.

Therefore Theorems \ref{thm:5} and \ref{thm:6} show that any
one-dimensional anisochronous mechanical system, when perturbed by a
periodic forcing and in the presence of dissipation, up to the
exceptional cases in which the functions $C_{k}(t_{0})$
are constant --- hence vanish, by Lemma \ref{lem:5} --- in $t_{0}$
for all $k\in\ZZZ_{+}$, admits subharmonic solutions of all orders,
without any assumption on the perturbation, --- a result which
does not follow from the analysis of \cite{HT1,CH}.

The case that all the functions $C_{k}(t_{0})$ are identically
constant in $t_{0}$ is really exceptional. This can be
appreciated by the following argument. If the function
$C(\eps,t_{0})$ does not depend on $t_{0}$ then not only,
by Lemma \ref{lem:5}, it must vanish identically,
i.e. $C(\eps,t_{0})=C(\eps)\equiv 0$, but we find also
that $t_{0}$ is left undetermined. In other words
the periodic solution persists for all values of $t_{0}$.
This means that if we take the system (\ref{eq:5.1}) with $\g=0$,
so that it becomes an autonomous quasi-integrable Hamiltonian system,
with no dissipation left, the full resonant torus with
frequency $\om=p/q$ persists under perturbation. This situation is
certainly unlikely --- even if not impossible in principle.
For instance one can take the system described by the Hamiltonian
\begin{equation}
H(x,y,t)=\frac{1}{2}y^{2} + \frac{1}{4} x^{4} +
\eps f(t) \left( \frac{1}{2}y^{2} + \frac{1}{4} x^{4} - E \right)^{2} ,
\label{eq:5.18} \end{equation}
with $E$ corresponding to the unperturbed solution $(x_{0}(t),y_{0}(t))$
with frequency $\om$. Then such a solution still satisfies the
corresponding Hamilton equations for all values of $\eps$
and for all values of the initial phase $t_{0}$: that is the full
resonant torus with frequency $\om$ persists. In particular
if $\om=p/q$ is rational --- so that the frequency of the
unperturbed solution becomes commensurate with the
frequency $1$ of the perturbing potential $f$,
the corresponding torus is resonant.

It is important to stress that if we look for a subharmonic solution
which continues some unperturbed periodic solution with a given
period $T=2\pi q/p$ it can happen that the corresponding
integral $\langle y_{0}f(x_{0}(\cdot), \cdot + t_{0}) \rangle$
identically vanishes. In fact, if $f$ is a trigonometric
polynomial (which is often the case in physical applications)
this happens for all $p/q$ but a finite set of values.
An explicit example has been considered in \cite{BBDGG}.
In these cases the subharmonic Melnikov function does not depend
on $t_{0}$ and it is linear in $C$: hence (\ref{eq:5.8})
can be satisfied only by taking $C_{0}(t_{0}) \equiv 0$.
Then, it becomes essential to go to higher
orders of perturbation theory to study for which values
of $C$ a subharmonic solution of order $q/p$ appears.
Again, we refer to \cite{BBDGG} for a situation in which one must
perform a higher order analysis to explain the numerical findings.

%%%%%%%%%%%%%%%%%%%%%%%%%%%%%%%%%%%%%%%%%%%%%%%%%%%%%%%%%%%%%%%%%%%%%%%%%
%%%%%%%%%%%%%%%%%%%%%%%%%%%%%%%%%%%%%%%%%%%%%%%%%%%%%%%%%%%%%%%%%%%%%%%%%
\zerarcounters
\section{Conclusions, and final comments}
\label{sec:6}
%%%%%%%%%%%%%%%%%%%%%%%%%%%%%%%%%%%%%%%%%%%%%%%%%%%%%%%%%%%%%%%%%%%%%%%%%
%%%%%%%%%%%%%%%%%%%%%%%%%%%%%%%%%%%%%%%%%%%%%%%%%%%%%%%%%%%%%%%%%%%%%%%%%

The Melnikov theory \cite{GH} considers systems which,
in suitable coordinates, can be written as in (\ref{eq:2.1}),
without the parameter $C$:
\begin{equation}
\begin{cases}
\dot \al = \om(A) + \eps F(\al,A,t) , \\
\dot A = \eps G(\al,A,t) , \end{cases}
\label{eq:6.1} \end{equation}
where all notations are as explained after (\ref{eq:2.1}).
Define the subharmonic Melnikov function as
\begin{equation}
M(t_{0}) = \frac{1}{T} \int_{0}^{T} {\rm d} t \,
G(\al_{0}(t),A_{0}(t),t+t_{0}) ,
\label{eq:6.2} \end{equation}
and set $M'(t_{0})={\rm d}M(t_{0})/{\rm d}t_{0}$.

We can repeat the analysis of formal solvability in Section \ref{sec:3},
with some adaptations due to the fact that no extra
parameters $C^{(k)}$ are at our disposal to any perturbation orders.

In particular to first order one needs $M(t_{0})=0$, so that
$t_{0}$ must be a zero for the subharmonic Melnikov function.
To higher orders we can write
\begin{equation}
G^{(k)}(\al(t),A(t),C,t+t_{0}) =
\partial_{1} G(\al_{0}(t),A_{0},t+t_{0}) \, \Bar \al^{(k)}
+ \Gamma^{(k)}(\al(t),A(t),t+t_{0}) ,
\label{eq:6.3} \end{equation}
where the function $\Gamma^{(k)}(\al(t),A(t),t+t_{0})$
depends on the corrections $\Bar \al^{(k')}$ to the initial phase,
only with $k'<k$.

To any perturbation order $k$ the constant $\Bar \al^{(k)}$ is
left undetermined. Anyway we are no longer free to fix it equal
to some arbitrary value, for instance zero, as we no longer have
the initial phase $t_{0}$ and the constants $C^{(k)}$ as free parameters.
Hence we shall need the corrections $\Bar \al^{(k)}$ to assure
solvability of the equations of motion to any order.
This will be possible in the light of the following results.

%%%%%%%%%%%%%%%%%%%%%%%%%%%%%%%%%%%%%%%%%%%%%%%%%%%%%%%%%%%%%%%%%%%%%%%%%
\begin{lemma} \label{lem:6}
One has $\om(A_{0}) \langle \partial_{1} G(\al_{0}(\cdot),
A_{0},\cdot+t_{0}) \rangle = - M'(t_{0})$.
\end{lemma}
%%%%%%%%%%%%%%%%%%%%%%%%%%%%%%%%%%%%%%%%%%%%%%%%%%%%%%%%%%%%%%%%%%%%%%%%%

\prova One has
\begin{equation}
\frac{{\rm d}}{{\rm d} t} G(\al_{0}(t),A_{0},t+t_{0}) =
\om(A_{0}) \, \partial_{1}
G(\al_{0}(t),A_{0},t+t_{0}) +
\frac{\partial}{\partial t_{0}} 
G(\al_{0}(t),A_{0},t+t_{0})  ,
\label{eq:6.4} \end{equation}
where we have used the fact that $\dot A_{0}(t)=0$ and
$\dot\al_{0}(t)=\om(A_{0})$. If we integrate (\ref{eq:6.4})
over a period we obtain
\begin{equation}
0 = \frac{1}{T} \int_{0}^{T} {\rm d} t \,
\frac{{\rm d}}{{\rm d} t} G(\al_{0}(t),A_{0},t+t_{0}) =
\om(A_{0}) \langle \partial_{1} G(\al_{0}(\cdot),A_{0},\cdot+t_{0})
\rangle + \frac{\partial}{\partial t_{0}}
\langle G(\al_{0}(\cdot),A_{0},\cdot+t_{0}) \rangle ,
\label{eq:6.5} \end{equation}
so that
\begin{equation}
\om(A_{0}) \partial_{1} G(\al_{0}(t),A_{0},C_{0},t+t_{0})
= - \frac{\partial}{\partial t_{0}}
\langle G(\al_{0}(t),A_{0},t+t_{0}) \rangle = - M'(t_{0}) .
\label{eq:6.6} \end{equation}
Hence the assertion follows. \qed

\*

Thus, if we impose the condition that $t_{0}$ be a simple zero for
the subharmonic Melnikov function we find that in (\ref{eq:6.3})
the derivative $\partial_{1} G(\al_{0}(t),A_{0},t+t_{0})$
is different from zero, and this allows us to fix $\Bar \al^{(k)}$
in such a way as to make the mean of $G^{(k)}(\al(t),A(t),t+t_{0})$
vanish. Hence by fixing the constants $\Bar A^{(k)}$ as explained
in Section \ref{sec:3} and the constants $\Bar \al^{(k)}$
as stated above we find that a solution in the form
of a formal power series in $\eps$ exists. The convergence
of the series, hence the existence of an analytic
solution, can be proved by reasoning as in Section \ref{sec:4}.
We do not repeat the analysis, which would essentially be
a word for word copy of what was done in Section \ref{sec:4}.

Therefore we have proved the following result --- well-known
in the literature \cite{GH}.

%%%%%%%%%%%%%%%%%%%%%%%%%%%%%%%%%%%%%%%%%%%%%%%%%%%%%%%%%%%%%%%%%%%%%%%%%
\begin{theorem} \label{thm:7}
Consider a periodic solution with frequency $\om=p/q$ for the
system (\ref{eq:6.1}), and assume that $t_{0}$ is a simple zero
for the subharmonic Melnikov function (\ref{eq:6.2}) corresponding to
such a solution. There exists $\eps_{0}>0$ such that for $|\eps|<\eps_{0}$
the system (\ref{eq:6.1}) has at least one subharmonic solution
of order $q/p$.
\end{theorem}
%%%%%%%%%%%%%%%%%%%%%%%%%%%%%%%%%%%%%%%%%%%%%%%%%%%%%%%%%%%%%%%%%%%%%%%%%

However, our analysis permits us to generalise  the result above. Define
\begin{equation}
M_{0}(t_{0})=M(t_{0}) , \qquad M_{k}(t_{0}) =
\langle \Gamma^{(k)}(\al(\cdot),A(\cdot),\cdot+t_{0}) \rangle ,
\quad k \in \NNN ,
\label{eq:6.7} \end{equation}
where the notations of (\ref{eq:6.3}) have been used. Note that if
$M_{k'}(t_{0})$ vanishes identically for all $k'=0,1,\ldots,k-1$,
then $M_{k}(t_{0})$ is well-defined. The following result follows.

%%%%%%%%%%%%%%%%%%%%%%%%%%%%%%%%%%%%%%%%%%%%%%%%%%%%%%%%%%%%%%%%%%%%%%%%%
\begin{theorem} \label{thm:8}
Consider a periodic solution with frequency $\om=p/q$ for the
system (\ref{eq:6.1}). Assume that the functions $M_{k'}$
are identically zero for all $k'=0,1,\ldots,k-1$, and assume that
$t_{0}$ is a simple zero for the function $M_{k}(t_{0})$.
There exists $\eps_{0}>0$ such that for $|\eps|<\eps_{0}$
the system (\ref{eq:6.1}) has at least one subharmonic solution
of order $q/p$.
\end{theorem}
%%%%%%%%%%%%%%%%%%%%%%%%%%%%%%%%%%%%%%%%%%%%%%%%%%%%%%%%%%%%%%%%%%%%%%%%%

Of course also the system (\ref{eq:2.1}) can be studied
as illustrated in this section. One simply treats the parameter
$C$ as fixed, and one fixes the initial phase $t_{0}$
in such a way that Theorem \ref{thm:7} or Theorem \ref{thm:8}
can be applied --- of course, provided the corresponding hypotheses
are satisfied. This has been done in \cite{BBDGG} to study
the subharmonic solutions of a forced cubic oscillator in
the presence of dissipation.

We also note that, as a particular case of Theorem \ref{thm:8},
it can happen that $M_{k}(t_{0}) \equiv 0$ for all $k\in\ZZZ_{+}$.
In that case formal solvability of the equations holds to all orders,
and the convergence of the series requires no condition on $t_{0}$,
and it can be proved by proceeding as in Section \ref{sec:4}.
In particular in such a case the full resonant torus persists
under perturbation. Of course, the vanishing identically of
all functions $M_{k}$ is a very unlikely situation, and, without
any further parameter at our disposal, we can hardly expect
this ever to happen. This shows that the persistence of the full torus
when the subharmonic Melnikov function is identically zero
is a very rare event.

The bifurcation curves studied in this paper concern subharmonic
solutions which are analytic in $\eps$. In principle
our results do not exclude existence of other subharmonic
solutions which are not analytic. Indeed, one could
wonder whether other periodic solutions with the same
period exist for $\eps\neq0$. In the presence of dissipation,
it is unlikely that solutions other than the attractive ones
found with the method we have
used, would be relevant for the dynamics
--- cf. for instance the problems
investigated in \cite{BDGG,BDG1,BDG2,BBDGG,BDGM}.
In general the situation can be delicate; for instance
when one investigates quasi-periodic solutions
corresponding to lower-dimensional tori of quasi-integrable
systems, where uniqueness becomes a subtle problem -- cf.
for instance \cite{GG2,CGGG}. Despite this, there are cases in
which the problem can be settled -- cf. \cite{BDG2,GCB}.

%%%%%%%%%%%%%%%%%%%%%%%%%%%%%%%%%%%%%%%%%%%%%%%%%%%%%%%%%%%%%%%%%%%%%%%%%%
%%%%%%%%%%%%%%%%%%%%%%%%%%%%%%%%%%%%%%%%%%%%%%%%%%%%%%%%%%%%%%%%%%%%%%%%%%
% References
%%%%%%%%%%%%%%%%%%%%%%%%%%%%%%%%%%%%%%%%%%%%%%%%%%%%%%%%%%%%%%%%%%%%%%%%%%
%%%%%%%%%%%%%%%%%%%%%%%%%%%%%%%%%%%%%%%%%%%%%%%%%%%%%%%%%%%%%%%%%%%%%%%%%%

\end{document}